\documentclass[a4paper,11pt,oneside]{article}

\usepackage[utf8]{inputenc}

\usepackage{ao-math-std}
\usepackage{ao-math-fields}

\usepackage{xcolor}
\ifcase1
\usepackage{luh-colors}
\colorlet{imgcol}{LUH-green}
\or
\colorlet{imgcol}{white}
\fi
\newcommand\image[2][]{\bgroup\fboxsep0pt\fboxrule0.2pt
  \protect\fcolorbox{imgcol}{imgcol!10}{\protect\includegraphics[#1]{#2}}\egroup}

\newcommand\qfor{\quad\forall}
\let\bs\boldsymbol
\newcommand\boldu{\bs u}
\newcommand\boldv{\bs v}
\newcommand\boldz{\bs z}
\newcommand\bolduu{\bs{uu}}
\newcommand\bolduz{\bs{uz}}
\newcommand\boldzu{\bs{zu}}
\newcommand\boldPhi{\bs\Phi}
\newcommand\bolddeltau{\bs{\delta u}}
\newcommand\bolddeltaz{\bs{\delta z}}
\newcommand\Phiu[1][]{\Phi_{u\optindex{#1}}}
\newcommand\Phiuy{\Phi_{u:y}}
\newcommand\Phiphi[1][]{\Phi_{\varphi\optindex{#1}}}
\newcommand\zu[1][]{z_{u\optindex{#1}}}

\newcommand\zphi[1][]{z_{\varphi\optindex{#1}}}
\newcommand\J{\mathcal J}
\newcommand\Lagr{\mathcal L}
\newcommand\dphi{\partial_t \varphi}
\newcommand\subphi[2]{_{\set{\varphi(#1) > \varphi(#2)}}}
\newcommand\subphim[2]{_{\set{\varphi^-_{#1} > \varphi^-_{#2}}}}
\newcommand\subdphi[1]{_{\set{\dphi(#1) > 0}}}
\newcommand\subdphiI[1]{_{\set{\dphi > 0, #1}}}
\newcommand\CIQ{W}
\newcommand\dtm[1][m]{\Delta t_#1}

\usepackage{amsmath,amsthm,amsfonts,amssymb,amscd}
\usepackage{mathtools}
\usepackage{booktabs}
\usepackage{pgfplots}
\pgfplotsset{compat=1.14,grid style={gray!15}}
\usepackage{siunitx}
\usepackage{float}

\newcommand\diff{\,\mathrm{d}}
\newcommand\dt{\diff t}
\newcommand\dx{\diff x}

\theoremstyle{plain} 
\newtheorem{theorem}{Theorem}[section]

\newtheorem{Remark}[theorem]{Remark}

\usepackage{authblk}

\usepackage[top=2cm, bottom=2cm, left=2cm, right=2cm]{geometry}

\usepackage{cleveref}

\begin{document}

\title{
Space-time formulation, discretization, and computations for
  phase-field fracture optimal control problems
}

\author[1]{D. Khimin}
\author[1]{M. C. Steinbach}
\author[1,2]{T. Wick}

\affil[1]{Leibniz Universit\"at Hannover, Institut f\"ur Angewandte
  Mathematik, Welfengarten 1, 30167 Hannover, Germany}

\affil[2]{Universit\'e Paris-Saclay, ENS Paris-Saclay, LMT -- Laboratoire de M\'ecanique et Technologie, 91190 Gif-sur-Yvette, France}

\date{}

\maketitle

\begin{abstract}
The purpose of this work is the development of space-time discretization
  schemes for phase-field optimal control problems.
  First, a time discretization of the forward problem is derived
  using a discontinuous Galerkin formulation.
  Here, a challenge is to include regularization terms
  and the crack irreversibility constraint.
  The optimal control setting is formulated by means of the Lagrangian approach
  from which the primal part, adjoint, tangent and adjoint Hessian are derived.
  Herein the overall Newton algorithm is based
  on a reduced approach by eliminating the state constraint.
  From the low-order discontinuous Galerkin discretization,
  adjoint time-stepping schemes are finally obtained.
  Our algorithmic developments are substantiated and illustrated
    with some numerical experiments.
\textbf{AMS}:\\
74R10, 65N30, 49M15, 49K20, 35Q74
\end{abstract}

\section{Introduction}

Fracture propagation using variational approaches and phase-field methods
is currently an important topic in applied mathematics and engineering.
The approach was established in \cite{FraMar98,BourFraMar00}
and overview articles and monographs include
\cite{BourFraMar08,BouFra19,WuNgNgSuBoSi19,Wi20_book,Fra21}
with numerous further references cited therein.
While the major amount of work concentrates on forward modeling
of phase-field fracture, more recently some work started
on parameter identification employing Bayesian inversion
\cite{KhoNoPaAbWiHei20,WuRoLoMa21,NoKhoUlAlWiFrWr21,NoKhoWi21},
stochastic phase-field modeling \cite{GERASIMOV2020113353},
and optimal control \cite{NeiWiWo17,NeiWiWo19,MoWo20}.

The main objective of this work is to design a computational framework
for the last topic mentioned,
namely phase-field fracture optimal control problems.
In prior work \cite{NeiWiWo17,NeiWiWo19} the emphasis was on mathematical analysis and
a brief illustration in terms of a numerical simulation for a fixed fracture.
However, computational details have not yet been discussed therein,
but are necessary in order to apply and investigate the methodology
for more practical applications such as propagating fractures.
Due to the irreversibility constraint on the fracture growth,
optimization problems subject to such an evolution become
mathematical programs with complementarity
constraints (MPCC) \cite{Barbu1984,MIGNOT1976,MIGNOT1984}
so that standard constraint qualifications like \cite{Robinson1976,Zowe1979} cannot hold.
Our computational approach requires stronger regularity
and hence we replace the complementarity constraint
with a suitable penalty term.

Designing a computational framework for phase-field fracture optimal
control is novel and challenging
because robust forward and optimization solvers are required.
For the forward solver, as intensively discussed in the literature, the linear
and nonlinear solutions are demanding because of the non-convexity of the
governing energy functional of the forward phase-field fracture model
and the relationship of discretization and regularization parameters.
For the nonlinear solution various methods were proposed
such as alternating minimization (staggered solution) \cite{Bour07,BuOrSue10},
quasi-monolithic solutions \cite{HeWheWi15,Wi20_book},
and fully monolithic schemes
\cite{GeLo16,Wi17_SISC,Wi17_CMAME,KoKr20,wambacq2020interiorpoint}.
Nonetheless, monolithic solutions remain difficult
and we add an additional viscous regularization term
as originally proposed in \cite{KneRoZa13}
and used in our governing model from \cite{NeiWiWo19}.
The optimization problem is formulated in terms of the reduced approach
by eliminating the state variable with a control-to-state operator.
Therein, Newton-type methods require the evaluation of the
\textit{adjoint, tangent}, and \textit{adjoint Hessian equations}.
The latter requires the evaluation of second-order derivatives;
see, e.g., \cite{BeMeVe07} for parabolic optimization problems.

The last paper serves as point of departure for our approach
in the current work.
Specifically, we employ Galerkin formulations in time
and discuss in detail how the crack irreversibility constraint is formulated
using a penalization \cite{MiWheWi13a,NeiWiWo17}
and an additional viscous regularization \cite{NeiWiWo19,KneRoZa13}.
Based on these settings, concrete time-stepping schemes are derived.
As usual, the primal and tangent problem run forward in time
whereas the adjoint and adjoint Hessian equations run backward in time.
We then adopt two numerical tests with a given initial notch (fracture)
in order to achieve a given fracture path while controlling Neumann
boundary traction forces. The main emphasis is to establish robust
numerical solver results in terms of the nonlinear forward solver
and the nonlinear optimization loop. We notice that propagating fractures
within numerical optimization are challenging and were not addressed
in the prior work \cite{NeiWiWo17,NeiWiWo19}.
Some further preliminary results
(yet with a stationary, non-propagating fracture) are published in
the book chapter \cite{KhiSteiWi21}.

The outline of this paper is as follows:
In \cref{sec_PFF}, the phase-field fracture forward model is introduced.
Furthermore, a Galerkin time discretization is provided.
Next, in \cref{sec_opt}, the optimization problem is stated,
including the reduced space approach.
In the key \cref{sec_Lag_aux} the Lagrangian and three auxiliary equations
are carefully derived in great detail. Then,
in \cref{sec_num_tests} two numerical experiments are discussed
in order to substantiate our algorithmic developments.
Our work is summarized in \cref{sec_conclusions}.

\section{Phase-field fracture forward model and space-time discretization}
\label{sec_PFF}
To formulate the forward problem, we first introduce some basic notation
and then proceed with a space-time discretization.

\subsection{Notation}
We consider a bounded domain $\Omega \subset \R^2$.
The boundary is partitioned as
$\partial \Omega = \Gamma_N \stackrel.\cup \Gamma_D$
where both $\Gamma_D$ and $\Gamma_N$ have nonzero Hausdorff measure.
Next we define two function spaces,
$V \coloneqq H_D^1(\Omega;\R^2) \times H^1(\Omega)$
for the displacement field $u$ and the phase-field $\varphi$,
and $Q \coloneqq L^2(\Gamma_N)$ for the control $q$, where
\begin{align*}
   H^1(\Omega;\R^2)
  &\coloneqq \defset{v \in L^2(\Omega;\R^2)}
    {D^{\alpha}v \in L^2(\Omega;\R^2) \ \forall
    \alpha \in \mathbb{N}_0^2, \ \abs\alpha \le 1}, \\
  H^1_D(\Omega; \R^2)
  &\coloneqq \defset{v \in H^1(\Omega; \R^2)}{v|_{\Gamma_D} = 0}.
\end{align*}
Moreover we consider a bounded time interval $I = [0, T]$
and introduce the spaces
\begin{equation*}
  X \coloneqq
  \defset{\boldu = (u, \varphi)}
  {\boldu \in L^2(I, V), \,
    \dphi \in L^2(I, H^{-1}(\Omega))}, \qquad
  \CIQ \coloneqq C(I, Q).
\end{equation*}
On $V$ respectively $X$ we define the scalar products
\begin{align*}
  (\boldu, \boldv)
  &\coloneqq \int_\Omega \boldu \cdot \boldv \dx
  \qfor \boldu, \boldv \in V, \\
  (\boldu, \boldv)_I
  &\coloneqq \int_I \int_\Omega \boldu \cdot \boldv \dx \dt
    = \int_I (\boldu(t), \boldv(t)) \dt
    \qfor \boldu, \boldv \in X,
\end{align*}
with induced norms $\norm{}$ and $\norm{}_I$,
and furthermore the restricted inner products
\begin{align*}
  (\boldu(t), \boldv(t))\subdphi{t}
  &\coloneqq \begin{cases}
    (\boldu(t), \boldv(t)), & \dphi(t) > 0, \\
    0, & \text{else},
  \end{cases} \\
  (\boldu, \boldv)\subdphiI{I}
  &\coloneqq \int_I (\boldu(t), \boldv(t))\subdphi{t} \dt
    \qfor \boldu, \boldv \in X,
\end{align*}
with induced semi-norms $\norm{}\subdphi{t}$
and $\norm{}\subdphiI{I}$.
We also notice that we later work with
$(\fcdot, \fcdot)\subphi{t_i}{t_j}$,
defined like $(\fcdot, \fcdot)\subdphi{t}$,
and with a semi-linear
form $a(\fcdot)(\fcdot)$ in which the first argument is nonlinear
and the second argument is linear.

\subsection{Energy functional of quasi-static variational fracture modeling}
In the next step we introduce a functional
$E_\varepsilon^\gamma \colon \CIQ \times X \to \R$
from which we derive our forward problem.
Here $E_\varepsilon^\gamma(q;u,\varphi)$ is defined as the sum
of the regularized total energy of a crack plus a penalty term
for the time dependent  irreversibility constraint $\dphi \le 0$.
The regularized total energy of a crack is given by
\begin{equation}\label{reg_total_energy}
  E_\varepsilon(q;u,\varphi) \coloneqq
  \frac12 (g(\varphi) \C e(u), e(u))_I -
  (q,u)_{\Gamma_N,I} + G_c \Gamma_\varepsilon(\varphi),
\end{equation}
where $q$ denotes a force that is applied in orthogonal direction to
$\Gamma_N \subset \partial \Omega$,
$\C$ is the elasticity tensor and $e(u)$ the symmetric gradient.
Then, we have
\[
  \C e(u) = \sigma(u) = 2 \mu e(u) + \lambda \operatorname{tr}(e(u)) I,
\]
where $\mu,\lambda>0$ are the Lam\'e parameters and $I$ is the identity matrix.
The so-called degradation function
$g(\varphi)\coloneqq (1-\kappa)\varphi^2 + \kappa$ helps to extend the displacements
to the entire domain $\Omega$.
The term $G_c\Gamma_\varepsilon(\varphi)\coloneqq \frac{1}{2\varepsilon} \norm{1-\varphi}_I+\frac{\varepsilon}{2} \norm{\nabla \varphi}_I^2 $
is a regularized form of the Hausdorff measure \cite{AmTo92}. So far the problem consists in finding a function $\boldu\coloneqq (u,\varphi) \in X$ that minimizes the regularized total energy $\eqref{reg_total_energy}$ subject to the irreversibility constraint $\dphi\leq 0$. In the sequel, the constraint is being replaced by a penalty term, which is defined as
\[
 R(\varphi)\coloneqq \norm{\dphi}\subdphiI{I}^2.
\]
In order to ensure differentiability up to second order, an alternative is to work
with a fourth-order penalization \cite{NeiWiWo17}.
One final modification of $E_\varepsilon$ is necessary.
We add the convexification term $\frac\eta2 \norm{\dphi}_I^2$
for some $\eta  > 0$.
Indeed, in \cite{NeiWiWo19}, the term $\eta (\varphi^i - \varphi^{i-1}, \psi)$
in time steps $i-1, i$ was considered for $\eta \ge 0$.
This term corresponds to a potential viscous regularization
of a rate-independent damage model \cite{KneRoZa13}.

Finally the forward problem consists in finding $\boldu = (u,\varphi) \in X$
that solves the following optimization problem for given intial data
$\boldu_0 = (u_0,\varphi_0) \in V$ and given control $q \in \CIQ$:
\begin{equation}
  \label{Cgamma}
  \min_{\boldu} \
  E_\varepsilon^{\gamma}(q;u,\varphi) \coloneqq
  E_\varepsilon(q;u,\varphi)
  + \frac{\gamma}{2} R(\varphi)
  + \frac\eta2 \norm{\dphi}_I^2,
\end{equation}
with penalty parameter $\gamma>0$ and convexification parameter $\eta > 0$.

\begin{Remark}[Initial condition $u_0$]
\label{rem_ini_u0}
Note that we are concerned with quasi-static brittle fracture
without explicit time derivative in the displacement equation.
Nonetheless, we introduce for formal reasons $u_0$.
First, we can develop in an analogous fashion time discretization schemes
for the overall forward model.
Second, it facilitates the extension to problems
in which the displacement equation does have a time derivative,
such as dynamic fracture \cite{BouLarRi11,BoVeScoHuLa12}.
Third, having $u_0$ allows for a monolithic implementation structure,
and the system matrix for the initial condition is regular.
\end{Remark}

\begin{Remark}[Convexification]
  We notice that strict positivity $\eta > 0$ ensures
  the required regularity in time, $\dphi \in L^2(I, H^{-1}(\Omega))$.
  Moreover, it improves the numerical solution process of \eqref{model}.
  In fact, one can show for the quasi-static case
  that for sufficiently large values of $\eta$ the control-to-state mapping
  associated with \eqref{Cgamma} is single valued
  due to strict convexity of the energy corresponding to the equation.
  However, the convexification term $\frac\eta2 \norm{\dphi}_I^2$
  also penalizes crack growth.
  To ensure the dominance of the physically motivated term
  $\frac{\gamma}{2} R(\varphi)$ we have to choose $\gamma \gg \eta$.
\end{Remark}

\subsection{Weak formulation}
Before we continue with the spatial discretization
and the concrete time-stepping scheme, we
state the weak form of \eqref{Cgamma}.
To this end we replace \eqref{Cgamma}
by its first order optimality conditions, see e.g., \cite{NeiWiWo17},
yielding a coupled nonlinear PDE system:
given $\boldu_0 \in V$ and $q \in \CIQ$,
find $\boldu \in X$ such that
\begin{equation}
  \label{model}
  \begin{aligned}
    (g(\varphi) \C e(u), e(\Phiu))_I - (q, \Phiu)_{\Gamma_N,I} &= 0, \\
    G_c \varepsilon (\nabla \varphi, \nabla \Phiphi)_I -
    \smash[b]{\frac{G_c}{\varepsilon}} (1 - \varphi, \Phiphi)_I +
    (1 - \kappa) (\varphi \C e(u) : e(u), \Phiphi)_I \\
    {} + \gamma (\dphi, \Phiphi)\subdphiI{I} + \eta (\dphi, \Phiphi)_I &= 0,
  \end{aligned}
\end{equation}
for every test function $\boldPhi = (\Phiu,\Phiphi) \in X$.

\subsection{Galerkin time discretization}
Using a time grid
$
  0 = t_0 < \dotsb < t_M = T,
$
we first partition the interval $I$ into $M$ left-open subintervals
$I_m = (t_{m-1}, t_m]$,
\[
  I = \set{0} \cup I_1 \cup \dotsb \cup I_M.
\]
By using the discontinuous Galerkin method, here dG(0), we seek for a solution $\boldu$
in the space $X^0_k$ of piecewise polynomials of degree $0$.
The subindex $k$ denotes the time-discretized function space
in order to distinguish from the continuous space $X$.
To this end, we have
\[
  X^0_k \coloneqq \defset{\boldv \in X}
  {\boldv|_{I_m} \in \polynomials_0(I_m, V), \,
    m = 1, \dots, M \text{ and } \boldv(0) \in V}.
\]

\begin{Remark}
\label{rem_constant_time}
Since we work with $r=0$, i.e., constant functions in time, we have
\[
  \partial_t \boldv = \boldv_m^- - \boldv_{m-1}^+ = 0
  \qfor \boldv \in X_k^0 \text{ and } m = 1, \dots, M.
\]
\end{Remark}

To work with the discontinuities in $X^0_k$, we introduce the notation
\[
  \boldv^+_m \coloneqq \boldv(t_m+), \qquad
  \boldv^-_m \coloneqq \boldv(t_m-) = \boldv(t_m), \qquad
  [\boldv]_m \coloneqq \boldv^+_m - \boldv^-_m.
\]
Now the discretized state equation can be derived from \eqref{model}
by combining the two equations into a single equation \eqref{State_dG}.
To simplify the notation let us replace the
energy-related terms of \eqref{model}
with a semi-linear form $a\colon Q \times V \times V \to \R$,
\begin{equation}
  \label{State_sum_4}
  \begin{aligned}
    a(q,\boldu)(\boldPhi)
    &\coloneqq g(\varphi) \cdot (\C e(u), e(\Phiu)) \\
    &\, + G_c \varepsilon (\nabla \varphi, \nabla \Phiphi)
    - \frac{G_c}{\varepsilon} (1 - \varphi, \Phiphi) \\
    &\, + (1 - \kappa) (\varphi \cdot \C e(u) : e(u), \Phiphi)
    - (q, \Phiuy)_{\Gamma_N}
    .
  \end{aligned}
\end{equation}
Here $\Phiuy$ denotes the $y$ component of $\Phiu =(\Phi_{u:x}, \Phiuy)$
in $\boldPhi = (\Phiu, \Phiphi) \in V$.
Now the fully discretized state equation consists of finding
a function $\boldu \in X^0_k$ for a given control $q$
such that for every $\boldPhi \in X^0_k$
\begin{subequations}
  \label{State_dG}
  \begin{align}
    \label{State_dG_a}
    0
    &= \sum_{m=1}^M \bigl[
      \gamma (\dphi, \Phiphi)\subdphiI{I_m}
      + \eta (\dphi, \Phiphi)_{I_m} \bigr] \\
    \label{State_dG_c}
    &+ \sum_{m=0}^{M-1} \bigl[
      \gamma ([\varphi]_m, \Phiphi[m]^+)\subphim{m+1}{m}
      + \eta ([\varphi]_m, \Phiphi[m]^+) \bigr] \\
    \label{State_dG_b}
    &+ \sum_{m=1}^M a(q(t_m),\boldu(t_m))(\boldPhi(t_m)) \dtm \\
    \label{State_dG_d}
    &+ (u^-_0 - u_0, \Phiu[0]^-) + (\varphi^-_0 - \varphi_0, \Phiphi[0]^-)
    .
  \end{align}
\end{subequations}
The time integral in \eqref{State_dG_b} has been approximated by the
right-sided box rule, where $\dtm \coloneqq t_m - t_{m-1}$.
Since the functions in $X^0_k$ might be discontinuous,
we have to add jump terms in the typical dG(0) manner,
which are contained in \eqref{State_dG_c}.
By expanding these jump terms,
\eqref{State_dG_c} (with index shifted by one) becomes
\begin{equation}
  \label{State_dG_tilde}
  \sum_{m=1}^M \bigl[
  \gamma (\varphi^+_{m-1} - \varphi^-_{m-1}, \Phiphi[m-1]^+)\subphim{m}{m-1}
  + \eta (\varphi^+_{m-1} - \varphi^-_{m-1}, \Phiphi[m-1]^+) \bigr].
\end{equation}
Now, since we are employing a dG(0) scheme, our test functions satisfy
\[
  \boldPhi^+_{m-1} = \boldPhi^-_m \qfor m = 1, \dots, M.
\]
Therefore \eqref{State_dG_a} vanishes entirely by \cref{rem_constant_time},
and the two terms containing $\varphi_{m-1}^+$
in \eqref{State_dG_tilde} can be written as
$(\varphi^-_m, \Phiphi[m]^-)\subphim{m}{m-1}$ and
$(\varphi^-_m, \Phiphi[m]^-)$, respectively.
Combining the resulting expression with
\eqref{State_dG_b} and \eqref{State_dG_d},
we finally rewrite \eqref{State_dG} as
\begin{equation}
  \label{State_sum_5}
  \begin{aligned}
    0
    &= \sum_{m=1}^M \Bigl( \gamma \bigl[
    (\varphi^-_m, \Phiphi[m]^-)\subphim{m}{m-1} -
    (\varphi^-_{m-1}, \Phiphi[m]^-)\subphim{m}{m-1} \bigr] \\[-2\jot]
    &\mkern70mu + \eta \bigl[
    (\varphi^-_m, \Phiphi[m]^-) -
    (\varphi^-_{m-1}, \Phiphi[m]^- ) \bigr] \\
    &\mkern70mu + a(q(t_m),\boldu(t_m))(\boldPhi(t_m)) \dtm \Bigr) \\
    &+ (u^-_0 - u_0, \Phiu[0]^-) + (\varphi^-_0 - \varphi_0, \Phiphi[0]^-)
    .
  \end{aligned}
\end{equation}

\subsection{Time-stepping scheme}
We begin the solution process by solving the last line of \eqref{State_sum_5}:
\begin{equation}
  \label{state_t=0}
  \begin{aligned}
    (u^-_0, \Phiu[0]^-) &= (u_0, \Phiu[0]^-), \\
    (\varphi^-_0, \Phiphi[0]^-) &= (\varphi_0, \Phiphi[0]^-),
  \end{aligned}
\end{equation}
or equivalently $(\boldu(0),\boldPhi_0) = (\boldu_0,\boldPhi_0)$.
Then we proceed and solve for $m = 1, \dots, M$ and every $\boldPhi \in X^0_k$
the following equation:
\begin{equation}
  \begin{aligned}\label{state_t>0}
    0
    &= \gamma (\varphi(t_m), \Phiphi(t_m))\subphi{t_m}{t_{m-1}}
    + \eta (\varphi(t_m), \Phiphi(t_m)) \\
    &- \gamma (\varphi(t_{m-1}), \Phiphi(t_{m}))\subphi{t_m}{t_{m-1}}
    - \eta (\varphi(t_{m-1}), \Phiphi(t_{m})) \\
    &+ a(q(t_m),\boldu(t_m))(\boldPhi(t_m)) \dtm
    .
  \end{aligned}
\end{equation}

\subsection{Spatial discretization}

For the spatial discretization, we employ again a Galerkin
finite element scheme by introducing $H^1$ conforming discrete spaces.
We consider two-dimensional shape-regular meshes
with quadrilateral elements $K$ forming  the mesh
$\mathcal{T}_h = \set{K}$; see \cite{Cia87}.
The spatial discretization parameter is the
diameter $h_K$ of the element $K$.
On the mesh $\mathcal{T}_h$ we construct a finite element space
$V_h \coloneqq V_{uh} \times V_{\varphi h}$ as usual:
\begin{align*}
  V_{uh} &\coloneqq
  \defset{v \in H^1_D(\Omega; \R^2)}
  {v|_K \in Q_s(K) \text{ for } K \in \mathcal{T}_h}, \\
  V_{\varphi h} &\coloneqq
  \defset{v \in H^1(\Omega)}
  {v|_K \in Q_s(K) \text{ for } K \in \mathcal{T}_h}.
\end{align*}
Herein $Q_s(K)$ consists of shape functions
that are obtained as bilinear transformations of functions
defined on the master element $\hat{K} = (0, 1)^2$,
where $\hat Q_s(\hat K)$ is the space
of tensor product polynomials up to degree $s$ in dimension $d$ defined as
\[
  \hat Q_s (\hat K) \coloneqq \operatorname{span}
  \Defset{\prod_{i=1}^d \hat x_i^{\alpha_i}}{\alpha_i \in \set{0, 1 \dots, s}}.
\]
Specifically, for $s = 1$ and $d = 2$ we have
\[
  \hat Q_1(\hat K)
  =
  \operatorname{span}\set{1, \hat x_1, \hat x_2, \hat x_1 \hat x_2}.
\]
We notice for the next two sections that the following derivations are independent
of the specific spatial discretization
and for this reason the subindex $h$ is omitted.

\section{Optimization with phase-field fracture}
\label{sec_opt}
We formulate the following separable NLP
with a tracking type cost functional.
For given $(u_0,\varphi_0) \in V$ we seek a solution
$(q,\boldu) \in \CIQ \times X_k^0$ of
\begin{equation}
\label{EQ:NLPgamma}
\begin{aligned}
  \min_{q,\textbf{u}} \quad & \J(q,\textbf{u}) \coloneqq
  \frac{1}{2}\sum_{m=1}^M  \norm{\varphi(t_m) - \varphi_d(t_m)}^2
  + \frac{\alpha}{2}\sum_{m=1}^M  \norm{q(t_m) - q_d(t_m)}_{\Gamma_N}^2
  \\
  \text{s.t.} \quad &(q,\boldu)
  \text{ solves \eqref{state_t=0} and \eqref{state_t>0} for  } m=1,\dots,M, 
\end{aligned}
\end{equation}
where $\varphi_d \in L^{\infty}(\Omega)$ is some desired phase-field
and $q_d$ is a suitable nominal control
  that we use for numerical stabilization.
The second sum represents a common Tikhonov regularization with the
Tikhonov parameter $\alpha$.
The existence of a global solution of $\eqref{EQ:NLPgamma}$
in $L^2(I,Q) \times X$
has been shown in \cite[Theorem 4.3]{NeiWiWo17}
for functions that are non-negative and weakly semi-continuous.

\subsection{Reduced optimization problem and solution algorithm}
We solve \eqref{EQ:NLPgamma} by a reduced space approach.
To this end, we assume the existence of a solution operator
$S\colon \CIQ \to X$ via equation \eqref{model}.
With this solution operator the cost functional $\J(q, \boldu)$
can be reduced to $j\colon \CIQ \to \R$, $j(q) \coloneqq \J(q, S(q))$.
As a result we can replace \eqref{EQ:NLPgamma} by the
unconstrained optimization problem
\begin{equation}
  \label{EQ:NLPgammared}
  \min_q \ j(q).
\end{equation}
The reduced problem is solved by Newton's method applied to $j'(q) = 0$,
and hence we need computable representations of the derivatives $j'$ and $j''$.
The established approach in \cite{BeMeVe07} requires
the solution of the following four equations
for the Lagrangian $\Lagr(q,\boldu,\boldz)$;
the concrete form is defined in \eqref{Lagrange}.
\begin{enumerate}
\item \emph{State equation:}
  given $q \in \CIQ$,
  find $\boldu \in X$ such that for all $\boldPhi \in X$
  \eqref{model} holds:
  \begin{equation}
    \label{Stateaux}
    \Lagr'_{\boldz}(q,\boldu,\boldz)(\boldPhi) = 0.
  \end{equation}
\item \emph{Adjoint equation:} given $q \in \CIQ$ and $\boldu = S(q)$,
  find $\boldz \in X$ such that for all $\boldPhi \in X$
  \begin{equation}
    \label{Adjointaux}
    \Lagr'_{\boldu}(q,\boldu,\boldz)(\boldPhi) = 0.
  \end{equation}
\item \emph{Tangent equation:} given $q \in \CIQ$, $\boldu = S(q)$
  and a direction $\delta q \in \CIQ$,
  find $\bolddeltau \in X$ such that for all $\boldPhi \in X$
  \begin{equation}
    \label{Tangentaux}
    \Lagr''_{q\boldz}(q,\boldu,\boldz)(\delta q, \boldPhi) +
    \Lagr''_{\boldu\boldz}(q,\boldu,\boldz)(\bolddeltau,\boldPhi) = 0.
  \end{equation}
\item \emph{Adjoint Hessian equation:} given $q \in \CIQ$,
  $\boldu = S(q)$, $\boldz \in X$ from \eqref{Adjointaux},
  $\bolddeltau \in X$ from \eqref{Tangentaux},
  and a direction $\delta q \in \CIQ$,
  find $\bolddeltaz \in X$ such that for all $\boldPhi \in X$
  \begin{equation}\label{AdjointHessianaux}
    \Lagr''_{q\boldu}(q,\boldu,\boldz)(\delta q,\boldPhi) +
    \Lagr''_{\bolduu}(q,\boldu,\boldz)(\bolddeltau,\boldPhi) +
    \Lagr''_{\boldzu}(q,\boldu,\boldz)(\bolddeltaz,\boldPhi) = 0.
  \end{equation}
\end{enumerate}
Solving these equations in a special order
(see for instance \cite{BeMeVe07,Me08})
leads to the following representations
of the derivatives that we need for Newton's method:
\begin{align*}
  j'(q)(\delta q)
  &= \Lagr'_q(q,\boldu,\boldz)(\delta q) \qfor \delta q \in \CIQ, \\
  j''(q)(\delta q_1, \delta q_2)
  &= \Lagr''_{qq}(q,\boldu,\boldz)(\delta q_1, \delta q_2) +
    \Lagr''_{\boldu q}(q,\boldu,\boldz)(\bolddeltau, \delta q_2) \\
  &+ \Lagr''_{\boldz q}(q,\boldu,\boldz)(\bolddeltaz,\delta q_2)
    \qfor \delta q_1, \delta q_2 \in \CIQ.
\end{align*}
\section{Lagrangian and auxiliary equations}
\label{sec_Lag_aux}
In the following main section, we specify the previously given
abstract formulations in detail. We first derive the Lagrangian
and then the three auxiliary equations
\eqref{Adjointaux}--\eqref{AdjointHessianaux}.
Specific emphasis is on the regularization terms
for the crack irreversibility and the convexification.

\subsection{Lagrangian}
We formulate the Lagrangian
$\Lagr\colon \CIQ \times X^0_k \times X^0_k \to \R$
within the dG(0) setting as
\begin{equation}
  \label{Lagrange}
  \begin{aligned}
    \Lagr(q,\boldu,\boldz)
    &\coloneqq \J(q,\boldu) \\
    &\,- \gamma (\dphi, \zphi)\subdphiI{I} - \eta (\dphi, \zphi)_I \\
    &\,- \int_I a(q(t),\boldu(t))(\boldz(t)) \dt \\
    &\,- \eta_0
    (u(0) - u_0, \zu(0)) -\eta (\varphi(0) - \varphi_0, \zphi(0))
    .
  \end{aligned}
\end{equation}
Note that we have scaled the initial conditions with
two different parameters $\eta_0$ and~$\eta$.
For the phase-field variable $\varphi$ we use the convexification parameter
of its time derivative to obtain $\eta (\varphi(0) - \varphi_0) = 0$.
This is common in the context of a dG(0) setting
as it produces desired cancelations with the jump terms
resulting from the discontinuities of the test functions.
In contrast, the initial condition for $u$ has no physical meaning.
Therefore we use a separate parameter $\eta_0 > 0$
to obtain $\eta_0 (u(0) - u_0) = 0$.
Later we choose $\eta_0 \ll \eta$.

\pagebreak[4]

\subsection{Adjoint}
In the adjoint for dG(0) we seek
$\boldz = (\zu, \zphi) \in X_k^0$ such that
\[
  \Lagr'_{\boldu}(q,\boldu,\boldz)(\boldPhi) = 0 \qfor \boldPhi \in X^0_k.
\]
The first interesting part is the calculation of the derivative of $\Lagr$.
We formulate it directly in the weak form
\begin{equation}
  \label{Lagrange_derivative_1}
  \begin{aligned}
    \Lagr'_{\boldu}(q,\boldu,\boldz) (\boldPhi)
    &= \J'_{\boldu}(q,\boldu)(\boldPhi) \\
    &- \gamma (\partial_t \Phiphi, \zphi)\subdphiI{I}
      - \eta (\partial_t \Phiphi, \zphi)_I \\
    &\,{}- \int_I a'_{\boldu}(q(t),\boldu(t))(\boldPhi(t),\boldz(t)) \dt \\
    &-\eta_0 (\Phiu(0), \zu(0)) -\eta (\Phiphi(0), \zphi(0))
    .
  \end{aligned}
\end{equation}
Herein the partial derivative of $a$ reads
\begin{equation}
  \label{a'_u}
  \begin{aligned}
    a'_{\boldu}(q,\boldu)(\boldPhi,\boldz)
    &= ((1-\kappa) \varphi^2 + \kappa) \cdot (\C e(\Phiu), e(\zu)) \\
    &+ 2 \varphi (1-\kappa) \Phiphi (\C e(u), e(\zu)) \\
    &+ G_c \varepsilon (\nabla \Phiphi, \nabla \zphi)
    + \frac{G_c}{\varepsilon} (\Phiphi, \zphi) \\
    &+ (1-\kappa) (\Phiphi \cdot \C e(u) : e(u), \zphi) \\
    &+ 2 \varphi (1-\kappa) (\C e(\Phiu) : e(u), \zphi)
    .
  \end{aligned}
\end{equation}
Now the main problem is that the time derivatives are applied
to the test function~$\boldPhi$ as usual in the adjoint.
Therefore we use integration by parts
to shift the time derivatives over to $\boldz$.
Then the second line
in \eqref{Lagrange_derivative_1} becomes
\begin{equation}
  \label{Lagrange_derivative_2}
  \begin{aligned}
    &\quad\gamma (\Phiphi, \partial_t \zphi)\subdphiI{I}
    + \eta (\Phiphi, \partial_t \zphi)_I \\
    &+ \gamma (\Phiphi(0), \zphi(0))\subdphi{0} + \eta (\Phiphi(0), \zphi(0)) \\
    &- \gamma (\Phiphi(T), \zphi(T))\subdphi{T} - \eta (\Phiphi(T), \zphi(T)) .
  \end{aligned}
\end{equation}
At this point we have to decide how to approximate
the time derivative $\dphi(0)$.
While $\dphi(t_m)$ for $m = 1, \dots, M$ is easily approximated
by the backward difference
\[
  \dphi(t_m) \approx \frac{\varphi(t_m) - \varphi(t_{m-1})}{t_m - t_{m-1}},
\]
this procedure will not work for the first mesh point $t_0 = 0$.
The forward difference
\[
  \dphi(0) \approx \frac{\varphi(t_1) - \varphi(t_0)}{t_1 - t_0}
\]
is a good choice because it simplifies the condition $\dphi(t_0) > 0$
to $\varphi(t_1) > \varphi(t_0)$ and leads to desired
cancelations in \eqref{Lagrange_derivative_dG}.
Now we will repeat the procedure that we applied to the state equation.
We approximate the time derivatives and add the jump terms
(with shifted index) as we did in \eqref{State_dG},
obtaining expressions similar to \eqref{State_dG_tilde}:
\begin{equation}
  \label{Lagrange_derivative_dG}
  \begin{aligned}
    \Lagr'_{\boldu}(q,\boldu,\boldz) (\boldPhi)
    &= \J'_{\boldu}(q,\boldu)(\boldPhi) \\
    &+ \sum_{m=1}^M \bigl[ \gamma
      (\Phiphi[m]^-, \zphi[m]^- - \zphi[m-1]^+)\subphim{m}{m-1} \\[-2\jot]
    &\mkern60mu + \eta (\Phiphi[m]^-, \zphi[m]^- - \zphi[m-1]^+) \bigr] \\
    &- \gamma (\Phiphi[M]^-, \zphi[M]^-)\subphi{t_M}{t_{M-1}} -
    \eta (\Phiphi[M]^-, \zphi[M]^-) \\
    &+ \gamma (\Phiphi[0]^-, \zphi[0]^-)\subphi{t_1}{t_0} +
    \eta (\Phiphi[0]^-, \zphi[0]^-) \\
    &+ \sum_{m=1}^M \bigl[ \gamma
      (\Phiphi[m-1]^-, \zphi[m-1]^+ - \zphi[m-1]^-)\subphim{m}{m-1} \\[-2\jot]
    &\mkern60mu + \eta (\Phiphi[m-1]^-, \zphi[m-1]^+ - \zphi[m-1]^-) \bigr] \\
    &- \sum_{m=1}^M
    a'_{\boldu}(q(t_m),\boldu(t_m))(\boldPhi(t_m),\boldz(t_m)) \dtm \\
    &-\eta_0 (\Phiu[0]^-, \zu[0]^-) -\eta (\Phiphi[0]^-, \zphi[0]^-)
    .
  \end{aligned}
\end{equation}
Since $\zphi\in X^0_k$, we have $\zphi[m]^- = \zphi[m-1]^+$
and see that the first sum vanishes entirely.
We also see that the terms $\pm \eta (\Phiphi[0]^-, \zphi[0]^-)$
in the fifth and the last line of \eqref{Lagrange_derivative_dG} cancel.
Moreover, we assume that $\varphi(t_1) \leq \varphi(t_0)$ in the initial step,
and hence the term $-\gamma (\Phiphi[0]^-, \zphi[0]^-)\subphi{t_1}{t_0}$
in the fifth line vanishes as well.

\begin{Remark}[Projection of the initial solution] \label{Remark_Projection_of_initial_solution}
  The assumption $\varphi(t_1) \leq \varphi(t_0)$ is numerically
  justified since at $t_0$ some initial phase-field solution is prescribed.
  From $t_0$ to $t_1$ an $L^2$ projection of the initial conditions
  is employed that conserves the crack irreversibility constraint.
\end{Remark}

By the above arguments we eliminate the second, third and fifth line
of \eqref{Lagrange_derivative_dG} and the second term of the last line,
whereas the initial values for $\zu$ are still present:
\begin{equation}
  \label{eq:Lu-final}
  \begin{aligned}
    \Lagr'_{\boldu}(q,\boldu,\boldz)(\boldPhi)
    &= \J'_{\boldu}(q,\boldu)(\boldPhi) \\
    &- \gamma (\Phiphi[M]^-, \zphi[M]^-)\subphi{t_M}{t_{M-1}}
      - \eta (\Phiphi[M]^-, \zphi[M]^-) \\
    &+ \sum_{m=1}^M \bigl[ \gamma
    (\Phiphi[m-1]^-, \zphi[m-1]^+ - \zphi[m-1]^-)\subphim{m}{m-1} \\[-2\jot]
    &\mkern60mu + \eta (\Phiphi[m-1]^-, \zphi[m-1]^+ - \zphi[m-1]^-) \bigr] \\
    &- \sum_{m=1}^M
      a'_{\boldu}(q(t_m),\boldu(t_m))(\boldPhi(t_m),\boldz(t_m)) \dtm \\
    &- \eta_0 (\Phiu[0]^-, \zu[0]^-)
    .
  \end{aligned}
\end{equation}

\newpage
\subsection{Adjoint time-stepping scheme}

From here on we exploit the separable structure of
$\J(q,\boldu) = \sum_m J(q(t_m),\boldu(t_m))$.
We start the solution process by pulling out from \eqref{eq:Lu-final}
every term associated with the last time point $t_M$:
\begin{equation}
  \begin{aligned}
    a_{\boldu}'&(q(t_M)\boldu(t_M))(\boldPhi(t_M),\boldz(t_M)) \dtm[M] \\
    &+ \gamma (\Phiphi[M]^-, \zphi[M]^-)\subphim{m}{m-1}
    + \eta (\Phiphi[M]^-, \zphi[M]^-) \\
    &= J_{\boldu}'(q(t_M),\boldu(t_M))(\boldPhi(t_M))
    \qfor \boldPhi \in X^0_k
    .
  \end{aligned}
\end{equation}
Now we collect what is left,
multiply by $-1$ and use the $X_k^0$ property ($\zphi[m-1]^+ = \zphi[m]^-$):
\begin{equation}
  \begin{aligned}
    0
    &= \sum_{m=1}^M \bigl[ \gamma
    (\Phiphi[m-1]^-, \zphi[m-1]^- - \zphi[m]^-)\subphim{m}{m-1}
    + \eta (\Phiphi[m-1]^-, \zphi[m-1]^- - \zphi[m]^-) \bigr] \\
    &+ \sum_{m=1}^{M-1}
    a'_{\boldu}(q(t_m),\boldu(t_m))(\boldPhi(t_m),\boldz(t_m)) \dtm \\
    &- \sum_{m=1}^{M-1}
    J'_{\boldu}(q(t_m),\boldu(t_m))(\boldPhi(t_m)) \\
    &+ \eta_0 (\Phiu[0]^-, \zu[0]^-)
    \qfor \boldPhi \in X_k^0.
 \end{aligned}
\end{equation}
To formulate the equations that are actually solved in every time step
we want to rewrite the entire equation as a single sum.
Therefore we shift down the index of the first sum (the jump terms),
take out the terms for $m = 0$, and obtain
\begin{equation*}
  \begin{aligned}
    0
    &=\sum_{m=1}^{M-1} \Bigl( \bigl[
    \gamma (\Phiphi[m]^-, \zphi[m]^- - \zphi[m+1]^-)\subphim{m+1}{m}
    + \eta (\Phiphi[m]^-, \zphi[m]^- - \zphi[m+1]^-) \bigr] \\[-2\jot]
    &\mkern70mu +
    a'_{\boldu}(q(t_m),\boldu(t_m))(\boldPhi(t_m),\boldz(t_m))
    \dtm \notag \\[\jot]
    &\mkern70mu -
    J'_{\boldu}(q(t_m),\boldu(t_m))(\boldPhi(t_m)) \Bigr) \\
    &+ \gamma (\Phiphi[0]^-, \zphi[0]^- - \zphi[1]^-)\subphim{1}{0}
    + \eta (\Phiphi[0]^-, \zphi[0]^- - \zphi[1]^-) \\
    &+ \eta_0 (\Phiu[0]^-, \zu[0]^-)
    .
  \end{aligned}
\end{equation*}
Now we solve for $m = M-1, M-2, \dots, 1$ the equation
\begin{align*}
  a_{\boldu}'&(q(t_m),\boldu(t_m))(\boldPhi(t_m),\boldz(t_m)) \dtm \\
  &+ \gamma (\Phiphi[m]^-, \zphi[m]^- - \zphi[m+1]^-)\subphim{m+1}{m}
    + \eta (\Phiphi[m]^-, \zphi[m]^- - \zphi[m+1]^-) \\
  &= J_{\boldu}'(q(t_m),\boldu(t_m))(\boldPhi(t_m))
  \qfor \boldPhi \in X^0_k
  .
\end{align*}
Finally three terms are left for $m = 0$,
\begin{equation}
  \label{Adjoint_t_0}
  \gamma (\Phiphi[0]^-, \zphi[0]^- - \zphi[1]^-)\subphim{1}{0}
  + \eta (\Phiphi[0]^-, \zphi[0]^- - \zphi[1]^-)
  + \eta_0 (\Phiu[0]^-,\zu[0]^-) = 0
  .
\end{equation}
For $\eta_0 \ll \eta$ small enough the last term of \eqref{Adjoint_t_0}
can be dropped and the following equation can be solved instead:
\begin{equation}\label{Adj_initial}
  (\Phiphi[0]^-, \zphi[1]^-) = (\Phiphi[0]^-, \zphi[0]^-).
\end{equation}

\begin{Remark}[Algorithmic realization]\label{Remark_Algorithmic_realization}
  To avoid singular matrices that would lead to a loss of convergence
  in the linear solvers,
  we have to add an intial condition for $\zu[0]^-$:
  $(\Phiu[0]^-, \zu[1]^-) = (\Phiu[0]^-, \zu[0]^-)$.
  In total we replace \eqref{Adj_initial} by
  $(\boldPhi_0^-, \boldz^-_1) =  (\boldPhi_0^-, \boldz^-_0)$.
  We also refer the reader to the third reason outlined in \cref{rem_ini_u0}.
\end{Remark}

\subsection{Tangent equation}
The second auxiliary equation is the tangent equation.
In this equation we seek
$\bolddeltau = (\delta u, \delta\varphi) \in X^0_k$ such that
\[
  \Lagr''_{q\boldz}(q,\boldu,\boldz)(\delta q,\boldPhi) +
  \Lagr''_{\bolduz}(q,\boldu,\boldz)(\bolddeltau,\boldPhi) = 0
  \qfor \boldPhi \in X_k^0.
\]
Here we will apply the same procedure as for the state equation.
Recall that $\Lagr(q,\boldu,\boldz)$ contains the integrand
$a(q(t),\boldu(t))(\boldz(t))$ with $\boldz(t)$ entering linearly.
Hence the partial derivative required for
$\Lagr''_{\bolduz}(q,\boldu,\boldz)(\bolddeltau,\boldPhi)$
is simply $a'_{\boldu}(q,\boldu)(\bolddeltau,\boldPhi)$,
and the partial derivative required for
$\Lagr''_{q\boldz}(q,\boldu,\boldz)(\delta q,\boldPhi)$
can be derived from \eqref{State_sum_4} as
\begin{equation}
  \label{a'_q}
  a'_q(q,\boldu)(\delta q, \boldPhi) = -(\delta q, \Phiuy)_{\Gamma_N}.
\end{equation}
Furthermore, $\J(q,\boldu)$ does not depend on $\boldz$,
hence $\J''_{q\boldz}$ and $\J''_{\boldu\boldz}$ vanish.
Using the right-sided box rule again,
we thus obtain the discretized tangent equation
\begin{equation}
  \label{Tangent}
  \begin{aligned}
    0
    &= \sum_{m=1}^M \bigl[ \gamma
    (\delta\varphi^-_m - \delta\varphi^+_{m-1}, \Phiphi[m]^-)\subphim{m}{m-1}
    + \eta (\delta\varphi^-_m - \delta\varphi^+_{m-1}, \Phiphi[m]^-) \bigr] \\
    &+ \sum_{m=1}^M
    a'_{\boldu}(q(t_m),\boldu(t_m))(\bolddeltau(t_m),\boldPhi(t_m)) \dtm \\
    &+ \sum_{m=0}^{M-1} \bigl[ \gamma
    (\delta\varphi^+_m - \delta\varphi^-_m, \Phiphi[m]^+)\subphim{m+1}{m}
    + \eta (\delta\varphi^+_m - \delta\varphi^-_m, \Phiphi[m]^+) \bigr] \\
    &{} + \eta_0 (\delta u^-_0, \Phiu[0]^-)
      + \eta(\delta\varphi^-_0, \Phiphi[0]^-) \\
    &+ \sum_{m=1}^M
    a'_q(q(t_m),\boldu(t_m))(\delta q(t_m), \boldPhi(t_m)) \dtm
    \qfor \boldPhi \in X^0_k
    .
  \end{aligned}
\end{equation}
It is clear that the first sum is zero due to the dG(0) property.
By shifting the index of the third sum in \eqref{Tangent}
and applying the dG(0) property to $\Phiphi[m-1]^+$
we can combine the last three sums and rewrite \eqref{Tangent} as
\begin{equation}
  \label{Tangent_reform}
  \begin{aligned}
    0
    &= \sum_{m=1}^M \Bigl(
    a'_{\boldu}(q(t_m),\boldu(t_m))(\bolddeltau(t_m),\boldPhi(t_m))
    \dtm \\[-2\jot]
    &\mkern70mu+ \gamma (\delta\varphi^+_{m-1} - \delta\varphi^-_{m-1},
    \Phiphi[m]^-)\subphim{m}{m-1} \\
    &\mkern70mu+
    \eta (\delta\varphi^+_{m-1} - \delta\varphi^-_{m-1}, \Phiphi[m]^-) \\
    &\mkern70mu+
    a'_q(q(t_m),\boldu(t_m))(\delta q(t_m), \boldPhi(t_m)) \dtm \Bigr) \\
    &{} + \eta_0 (\delta u^-_0, \Phiu[0]^-)
      + \eta(\delta\varphi^-_0, \Phiphi[0]^-)
    \qfor \boldPhi \in X^0_k
    .
  \end{aligned}
\end{equation}
\subsection{Tangent time-stepping schemes}

As in the state equation we first solve the initial conditions,
\[
\begin{aligned}
  (\delta u^-_0, \Phiu[0]^-) &= 0, \\
  (\delta\varphi^-_0, \Phiphi[0]^-) &= 0.
\end{aligned}
\]
Applying the $X_k^0$ property to $\delta \varphi_{m-1}^{+}$ we can finally solve for $m=1,\dots,M$
the following equation
\begin{equation}
  \label{Tanget_ste}
  \begin{aligned}
    \gamma(\delta\varphi^-_m, &\Phiphi[m]^-)\subphim{m}{m-1}
    + \eta (\delta\varphi^-_m, \Phiphi[m]^-) \\
    &+ a'_{\boldu}(q(t_m),\boldu(t_m))(\bolddeltau(t_m),\boldPhi(t_m)) \dtm \\
    &= (\delta\varphi^-_{m-1}, \Phiphi[m]^-) +
    (\delta\varphi^-_{m-1}, \Phiphi[m]^-)\subphim{m}{m-1} \\
    &- a'_q(q(t_m),\boldu(t_m))(\delta q(t_m),\boldPhi(t_m)) \dtm
    \qfor \boldPhi \in X^0_k
    .
  \end{aligned}
\end{equation}
\subsection{Adjoint Hessian equation}
The third and last auxiliary equation is the adjoint Hessian equation.
In this equation we seek $\bolddeltaz = (\delta\zu, \delta\zphi) \in X^0_k$
such that for all $\boldPhi \in X_k^0$ the following equation holds true:
\begin{equation}
  \label{Adj_Hessian}
  \Lagr''_{q\boldu}(q,\boldu,\boldz)(\delta q,\boldPhi) +
  \Lagr''_{\bolduu}(q,\boldu,\boldz)(\bolddeltau,\boldPhi) +
  \Lagr''_{\boldzu}(q,\boldu,\boldz)(\bolddeltaz, \boldPhi) = 0
  .
\end{equation}
First we see that $\Lagr''_{q\boldu}(q,\boldu,\boldz)(\delta q,\boldPhi) = 0$
since $q$ and $\boldu$ are decoupled. The derivative of $a$ in
$\Lagr''_{\boldzu}(q,\boldu,\boldz)(\bolddeltaz, \boldPhi)$
is given by $a'_{\boldu}(q,\boldu)(\boldPhi,\bolddeltaz)$
due to the linearity of $\boldz$ in $a$.
However, a genuine second-order derivative of $a$ arises in
$\Lagr''_{\bolduu}(q,\boldu,\boldz)(\bolddeltau,\boldPhi)$:
\begin{equation}
  \label{a'_uu}
  \begin{aligned}
    a''_{\bolduu}(q,\boldu)(\bolddeltau,\boldPhi,\boldz)
    &= 2 \varphi \cdot (1-\kappa) \Phiphi \cdot (\C e(\delta u), e(\zu)) \\
    &+ 2 \delta\varphi \cdot (1-\kappa) (\C e(u), e(\zu)) \cdot \Phiphi \\
    &+ 2 \varphi \cdot (1-\kappa) (\C e(u), e(\zu)) \delta\varphi \\
    &+ 2 \varphi \cdot (1-\kappa) (\C e(\Phiu) : e(\delta u), \zphi) \\
    &+ 2 \delta\varphi \cdot (1-\kappa) (\C e(\Phiu) : e(u), \zphi) \\
    &+ 2 (\C e(\delta u) : e(u), \zphi) \cdot \Phiphi
    .
  \end{aligned}
\end{equation}
Now we can rewrite \eqref{Adj_Hessian} in a dG(0) setting:
\begin{equation}
  \begin{aligned}
    0
    &= \sum_{m=1}^M
    J''_{\bolduu}(q(t_m),\boldu(t_m))(\bolddeltau(t_m),\boldPhi(t_m)) \\
    &- \sum_{m=1}^M a''_{\bolduu}(q(t_m),\boldu(t_m))
    (\bolddeltau(t_m),\boldPhi(t_m),\boldz(t_m)) \dtm \\
    &+ \sum_{m=1}^M \bigl[ \gamma
    (\Phiphi[m]^-, \delta\zphi[m]^- - \delta\zphi[m-1]^+)\subphim{m}{m-1}
    + \eta (\Phiphi[m]^-, \zphi[m]^- - \zphi[m-1]^+) \bigr] \\
    &- \gamma (\Phiphi[M]^-, \delta\zphi[M]^-)\subphim{M}{M-1}
    - \eta (\Phiphi[M]^-, \delta\zphi[M]^-) \\
    &+ \gamma (\Phiphi[0]^-, \delta\zphi[0]^-)\subphim{1}{0}
    + \eta(\Phiphi[0]^-, \delta\zphi[0]^-) \\
    &- \sum_{m=1}^M
    a'_{\boldu}(q(t_m),\boldu(t_m))(\boldPhi(t_m),\bolddeltaz(t_m)) \dtm \\
    &+\sum_{m=0}^{M-1} \gamma
    (\Phiphi[m]^-, \delta\zphi[m]^+ - \delta\zphi[m]^-)\subphim{m+1}{m}
    + \eta (\Phiphi[m]^-, \delta\zphi[m]^+ - \delta\zphi[m]^-) \\
    &- \eta_0 (\Phiu[0]^-, \delta\zu[0]^-)
    - \eta (\Phiphi[0]^-, \delta\zphi[0]^-)
    \qfor \boldPhi \in X^0_k
    .
  \end{aligned}
\end{equation}
Note that the same scaling of initial data was applied
that we already used for the adjoint equation.
By the $X_k^0$ property the third sum vanishes entirely.
Due to \cref{Remark_Projection_of_initial_solution} and the cancelation of
$\pm \eta(\Phiphi[0]^-, \delta\zphi[0]^-)$
the fifth line vanishes as well.
By shifting the index of the jump terms we can rewrite the equation as:
\begin{equation}
  \label{Adj_Hessian_1}
  \begin{aligned}
    0
    &= \sum_{m=1}^M \Bigl(
    J''_{\bolduu}(q(t_m),\boldu(t_m))(\bolddeltau(t_m),\boldPhi(t_m)) \\[-2\jot]
    &\mkern70mu- a''_{\bolduu}(q(t_m),\boldu(t_m))
    (\bolddeltau(t_m),\boldPhi(t_m),\boldz(t_m)) \dtm \\
    &\mkern70mu-
    a'_{\boldu}(q(t_m),\boldu(t_m))(\boldPhi(t_m),\bolddeltaz(t_m)) \dtm \\
    &\mkern70mu+ \gamma
    (\Phiphi[m-1]^-, \delta\zphi[m-1]^+ - \delta\zphi[m-1]^-)\subphim{m}{m-1} \\
    &\mkern70mu+ \eta
    (\Phiphi[m-1]^-, \delta\zphi[m-1]^+ - \delta\zphi[m-1]^-) \Bigr) \\
    &- \gamma (\Phiphi[M]^-, \delta\zphi[M]^-)\subphim{M}{M-1}
    - \eta(\Phiphi[M]^-, \delta\zphi[M]^-) \\
    &- \eta_0 (\Phiu[0]^-, \delta\zu[0]^-)
    \qfor \boldPhi \in X^0_k
    .
  \end{aligned}
\end{equation}

\subsection{{Adjoint Hessian time-stepping schemes}}

As in the adjoint time-step\-ping scheme
we first collect all terms that contain the last time point $t_M$ and solve
\begin{equation}
  \begin{aligned}
    0
    &= J''_{\bolduu}(q(t_M)\boldu(t_M))(\bolddeltau(t_M), \boldPhi(t_M)) \\
    &- a'_{\boldu}(q(t_M)\boldu(t_M))(\boldPhi(t_M),\bolddeltaz(t_M)) \dtm[M] \\
    &- a''_{\bolduu}(q(t_M)\boldu(t_M))
    (\bolddeltau(t_M),\boldPhi(t_M),\boldz(t_M)) \dtm[M] \\
    &- \gamma (\Phiphi[M]^-, \delta\zphi[M]^-)\subphim{M}{M-1}
    - \eta (\Phiphi[M]^-, \delta\zphi[M]^-)
    \qfor \boldPhi \in X^0_k
    .
  \end{aligned}
\end{equation}
Then \eqref{Adj_Hessian_1} becomes
\begin{equation}
  \label{Adj_Hessian_2}
  \begin{aligned}
    0
    &= \sum_{m=1}^{M-1} \Bigl(
    J''_{\bolduu}(q(t_m),\boldu(t_m))(\bolddeltau(t_m),\boldPhi(t_m)) \\[-2\jot]
    &\mkern70mu- a''_{\bolduu}(q(t_m),\boldu(t_m))
    (\bolddeltau(t_m),\boldPhi(t_m),\boldz(t_m)) \dtm \\
    &\mkern70mu- a'_{\boldu}(q(t_m),\boldu(t_m))(\boldPhi(t_m),\bolddeltaz(t_m))
    \dtm \Bigl) \\
    &+ \sum_{m=1}^M \Bigl( \gamma
    (\Phiphi[m-1]^-, \delta\zphi[m-1]^+ - \delta\zphi[m-1]^-)\subphim{m}{m-1}
    \\[-2\jot]
    &\mkern70mu+ \eta
    (\Phiphi[m-1]^-, \delta\zphi[m-1]^+ - \delta\zphi[m-1]^-) \Bigr) \\
    &- \eta_0 (\Phiu[0]^-, \delta\zu[0]^-)
    \qfor \boldPhi \in X^0_k
    .
  \end{aligned}
\end{equation}
In the final reformulation we shift the index of the second sum (jump-terms)
and take out the terms corresponding to $m = 0$
\begin{equation}
  \label{Adj_Hessian_3}
  \begin{aligned}
    0
    &= \sum_{m=1}^{M-1} \Bigl(
    J''_{\bolduu}(q(t_m),\boldu(t_m))(\bolddeltau(t_m),\boldPhi(t_m)) \\[-2\jot]
    &\mkern70mu- a''_{\bolduu}(q(t_m),\boldu(t_m))
    (\bolddeltau(t_m),\boldPhi(t_m),\boldz(t_m)) \dtm \\
    &\mkern70mu-
    a'_{\boldu}(q(t_m),\boldu(t_m))(\boldPhi(t_m),\bolddeltaz(t_m)) \dtm \\
    &\mkern70mu+ \gamma
    (\Phiphi[m]^-, \delta\zphi[m]^+ - \delta\zphi[m]^-)\subphim{m+1}{m} \\
    &\mkern70mu+ \eta
    (\Phiphi[m]^-, \delta\zphi[m]^+ - \delta\zphi[m]^-) \Bigr) \\
    &+ \gamma (\Phiphi[0]^-, \delta\zphi[0]^+ - \delta\zphi[0]^-)\subphim{1}{0}
    + \eta (\Phiphi[0]^-, \delta\zphi[0]^+ - \delta\zphi[0]^-) \\
    &- \eta_0 (\Phiu[0]^-, \delta\zu[0]^-)
    \qfor \boldPhi \in X^0_k
    .
  \end{aligned}
\end{equation}
As already pointed out in the time-stepping scheme for the adjoint equation,
all dual equations have to be solved backwards in time.
Therefore, we solve the following equation for $m=M-1,M-2,\dots,1$
\begin{equation*}
  \begin{aligned}
    0
    &= J''_{\bolduu}(q(t_m),\boldu(t_m))(\bolddeltau(t_m),\boldPhi(t_m)) \\
    &- a''_{\bolduu}(q(t_m),\boldu(t_m))
    (\bolddeltau(t_m),\boldPhi(t_m),\boldz(t_m)) \dtm \\
    &- a'_{\boldu}(q(t_m),\boldu(t_m))(\boldPhi(t_m),\bolddeltaz(t_m)) \dtm \\
    &+ \gamma
    (\Phiphi[m]^-, \delta\zphi[m]^+ - \delta\zphi[m]^-)\subphim{m+1}{m} \\
    &+ \eta (\Phiphi[m]^-, \delta\zphi[m]^+ - \delta\zphi[m]^-)
    \qfor \boldPhi \in X^0_k
    .
  \end{aligned}
\end{equation*}
As a result, the only remaning terms in \eqref{Adj_Hessian_3} are
\begin{equation}\label{Adj_hessian_before_drop}
  \gamma (\Phiphi[0]^-, \delta\zphi[0]^+ - \delta\zphi[0]^-)\subphim{1}{0}
  + \eta (\Phiphi[0]^-, \delta\zphi[0]^+ - \delta\zphi[0]^-)
  - \eta_0 (\Phiu[0]^-, \delta\zu[0]^-)
  .
\end{equation}
Finally we can apply the assumption $\eta_0 \ll \eta$ once more and drop the last term
in \eqref{Adj_hessian_before_drop}. Consequently the following equations have to be solved
for all $ \boldPhi \in X^0_k$:
\begin{align*}
  (\Phiphi[0]^-, \delta\zphi[0]^-) &= (\Phiphi[0]^-, \delta\zphi[1]^-), \\
  (\Phiu[0]^-, \delta\zu[0]^-) &= (\Phiu[0]^-, \delta\zu[1]^-).
\end{align*}
Note that \cref{Remark_Algorithmic_realization} was applied
to \eqref{Adj_hessian_before_drop} as well.

\newpage
\section{Numerical tests}
\label{sec_num_tests}
In the following section we present two numerical examples
for the optimal control problem.
In these examples we use the tracking type functional of \eqref{EQ:NLPgamma}
to find an optimal control that approximately produces a desired phase-field.
All numerical computations are performed with the open source software libraries
\texttt{deal.II} \cite{dealII91,deal2020}
and \textsc{DOpElib} \cite{dope,DOpElib}.

\newcommand\figex[1]{%
  \begin{tikzpicture}[scale=1]
    \draw (0,0)
    -- (0,4) node[pos=0.5,left]  {$\Gamma_{\text{free}}$}
    -- (4,4) node[pos=0.5,above] {$\Gamma_N$}
    -- (4,0) node[pos=0.5,right] {$\Gamma_{\text{free}}$}
    -- cycle node[pos=0.5,below] {$\Gamma_D$};
    \node[teal] at (2,1) {Example #1};
    \node at (2,3) {$\Omega$};
    \node[blue,above] at (3,4) {$\uparrow q$};
    \draw[dotted,red] (2-#1,2) -- (3-#1,2) node[pos=0.5,above,red] {$\varphi_d$};
    \draw[thick,blue] (3-#1,2) -- (5-#1,2) node[pos=0.5,above,blue] {notch};
    \fill
    (2-#1,2) circle (1.5pt) (3-#1,2) circle (1.5pt) (5-#1,2) circle (1.5pt);
  \end{tikzpicture}%
}

\begin{figure}[H]
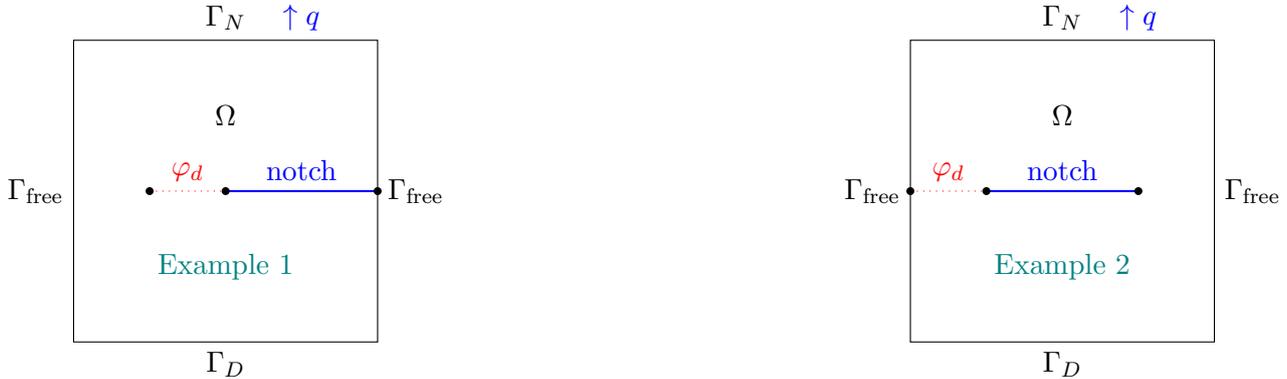

  \label{fig:examples}
  \centering
  \figex1\hfill\figex2%
  \caption{Domain $\Omega = (0, 1)^2$ with partitioned boundary
    $\partial \Omega$, initial notch, and desired crack $\varphi_d$.}
\end{figure}

For both examples we consider the square domain $\Omega = (0, 1)^2$
with a horizontal notch, see \cref{fig:examples}.
\begin{table}[H]
  \label{tab:params}
  \caption{Regularization and penalty parameters (left)
    and model parameters (right).}
  \centering
  \begin{tabular}{llS[table-format=1.2e+2]}
    \toprule
    Par. & Definition & {Value} \\
    \midrule
    $\varepsilon_1$ & Regul. (crack) $\approx 4 h_1$ & 0.0884 \\
    $\varepsilon_2$ & Regul. (crack) $\approx 4 h_2$ & 0.0442 \\
    $\kappa$ & Regul. (crack) & 1.00e-10 \\
    $\eta$ & Regul. (viscosity) & 1.00e3 \\
    $\gamma$ & Penalty & 1.00e5  \\
    $\alpha_1$ & Tikhonov & 4.75e-10 \\
    $\alpha_2$ & Tikhonov & 1.00e-10 \\
    \bottomrule
  \end{tabular}
  \hfill
  \begin{tabular}{llS[table-format=1.1e1]}
    \toprule
    Par. & Definition & {Value} \\
    \midrule
    $G_c$ & Fracture toughness & 1.0 \\
    $\nu_s$ & Poisson's ratio & 0.2 \\
    $E$ & Young's modulus & 1.0e6 \\ \\
    $q_0$ & Initial control & 1.0 \\
    $q_{d1}$ & Nominal control & 1.0e3 \\
    $q_{d2}$ & Nominal control & 3.0e3 \\
    \bottomrule
  \end{tabular}
\end{table}%
In Example~1 the notch is in the middle of the right side,
defined as $(0.5, 1) \times \set{0.5}$,
in Example~2 it is in the middle of $\Omega$,
defined as $(0.25, 0.75) \times \set{0.5}$.
The boundary $\partial \Omega$ is partitioned as
$\partial \Omega = \Gamma_N \cup \Gamma_D \cup \Gamma_{\text{free}}$,
where $\Gamma_N \coloneqq [0,1] \times \set{1}$,
$\Gamma_D \coloneqq [0, 1] \times \set{0}$, and
$\Gamma_{\text{free}} \coloneqq \set{0, 1} \times (0, 1)$.
On $\Gamma_N$ we apply the force $q$ in orthogonal direction to the domain
and on $\Gamma_D$ we enforce Dirichlet boundary conditions
for the displacement $u = 0$.
We choose the time interval $[0, 1]$
with 41 equidistant time points $t_m$,
i.e.\ $T = 1$ and $M = 40$.
The control space $Q_h$ (the spatial discretization of $Q$)
is one-dimensional in the sense that the force is only applied
in $y$-direction and is constant in time, $q(t_m) = q$.
The spatial mesh consists of $64 \times 64$ square elements in Example~1
and $128 \times 128$ square elements in Example~2,
hence the element diameter is
$h_1 = \sqrt2/64$ and $h_2 = \sqrt2/128$, respectively.
The initial data is given by $\boldu_0 = (u_0, \varphi_0)$
where $\varphi_0$ describes the horizontal notch,
\begin{equation}
  \label{EQ:phi0}
  \varphi_0(x, y) \coloneqq
  \begin{cases}
    0, & x \in (0.50, 1.00) \text{ and } y = 0.5 \text{ (Example 1)}, \\
    0, & x \in (0.25, 0.75) \text{ and } y = 0.5 \text{ (Example 2)}, \\
    1, & \text{else}.
  \end{cases}
\end{equation}
The desired phase-field $\varphi_d$ continues the initial notch to the left,
\begin{equation}
  \label{EQ:phid}
  \varphi_d(x, y) \coloneqq
  \begin{cases}
    0, & x \in (0.25, 0.5) \text{ and }
    y \in (0.5 - h_1, 0.5 + h_1) \text{ (Example 1)}, \\
    0, & x \in (0, 0.25) \text{ and }
    y \in (0.5 - 2 h_2, 0.5 + 2 h_2) \text{ (Example 2)}, \\
    1, & \text{else}.
  \end{cases}
\end{equation}
The parameters used in our numerical tests are given in \cref{tab:params}.
Note that $\varepsilon$, $\alpha$, and the constant nominal control $q_d$
differ for the two examples
while the crack widths that we prescribe via $\varphi_d$ agree:
$2 h_2 = h_1$ in \eqref{EQ:phid}.

\subsection{Example 1: horizontal fracture in right half domain}
\label{SEC:Example1}

The first example is motivated by a standard problem:
the single edge notched tension test \cite{MieWelHof10a,MieWelHof10b};
see again \cref{fig:examples}.
Our results are presented in \cref{tab:results1}.
The first column (Iter) gives the iteration index of Newton's method
in solving the reduced problem \eqref{EQ:NLPgammared}.
The second column (CG) gives the number of CG iterations required
for computing the Newton increment.
The remaining values are the relative and absolute Newton residuals,
the cost functional $\J$ and its tracking part
$\frac12 \sum_{m=1}^M \norm{\varphi(t_m) - \varphi_d(t_m)}^2$,
the maximal force $\abs{q_{\max}}$ applied on $\Gamma_N$,
and finally the Tikhonov regularization term,
$\frac{\alpha}{2}\sum_{m=1}^M  \norm{q(t_m) - q_d(t_m)}_{\Gamma_N}^2$.
All values are rounded to three or five significant digits.

\begin{figure}
  \label{fig:force1}
  \footnotesize
  \centering
  \begin{tikzpicture}
    \begin{axis}[width=0.8\linewidth,height=48mm,
      xmin=0,xmax=1,ymin=0,ymax=2500,ytick distance={500},grid=major]
      \draw[red,dotted] (0,1000) -- (1,1000) node[pos=0.48,below] {$q_d$};
      \addplot[blue] table {Figures/Test_1/Kraft.txt} node[pos=0.48,above] {$q$};
    \end{axis}
  \end{tikzpicture}
  \caption{Example 1: optimal control force (solid) and nominal control
    force (dotted) on upper boundary $\Gamma_N = [0, 1] \times \set{1}$.}
\end{figure}
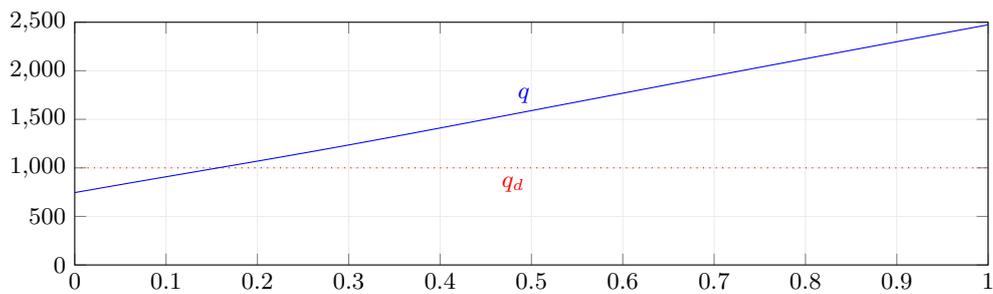

The Newton iteration terminates when the relative residual
or the absolute residual falls below the tolerance \num{2e-12}.
Optimal phase-fields, displacements, adjoints, and forces
are presented in \cref{fig:force1,fig:phase1,fig:displacement1,fig:adjoint1}.
The optimal force on $\Gamma_N = [0, 1] \times \set{1}$
grows almost linearly from 745.3 at $(0, 1)$ to 2473.4 at $(1, 1)$.

\begin{table}
  \label{tab:results1}
  \centering
  \caption{Results of Example 1: numerical solver performance,
      convergence of functionals, and evolution of control forces.}
  \sisetup{table-format=1.4e+1}
  \begin{tabular}{rr
    S[table-format=1.2e+1]S[table-format=1.2e+2]SSSS[table-format=4.1]}
    \toprule
    Iter&CG&{Relative}&{Absolute}&{Cost}&{Tracking}&{Tikhonov}&{Force}\\[-2pt]
        &  &{residual}&{residual}&      &          &       &          \\
    \midrule
    0 & ---& 1.0     & 4.62e-07 & 5.1630e-3 & 4.9289e-3 & 2.3406e-4 & 1.0\\
    1 &  3 & 0.464   & 2.14e-07 & 4.7699e-3 & 4.7699e-3 & 3.2317e-9 & 1001.2\\
    2 &  3 & 0.207   & 9.56e-08 & 4.5782e-3 & 4.5276e-3 & 5.0571e-5 & 1899.7\\
    3 &  3 & 0.106   & 4.88e-08 & 4.4968e-3 & 4.3911e-3 & 1.0570e-4 & 2308.1\\
    4 &  3 & 2.73e-3 & 1.26e-09 & 4.4600e-3 & 4.3186e-3 & 1.4137e-4 & 2498.9\\
    5 & 10 & 7.18e-4 & 3.32e-10 & 4.4606e-3 & 4.3203e-3 & 1.4032e-4 & 2469.5\\
    6 &  5 & 2.93e-4 & 1.35e-10 & 4.4592e-3 & 4.3174e-3 & 1.4183e-4 & 2471.1\\
    7 &  5 & 1.03e-4 & 4.73e-11 & 4.4585e-3 & 4.3159e-3 & 1.4258e-4 & 2473.1\\
    8 &  3 & 4.44e-5 & 2.05e-11 & 4.4584e-3 & 4.3156e-3 & 1.4273e-4 & 2473.0\\
    9 &  3 & 1.54e-5 & 7.09e-12 & 4.4583e-3 & 4.3154e-3 & 1.4286e-4 & 2473.4\\
    10 & 2 & 3.78e-6 & 1.75e-12 & 4.4582e-3 & 4.3153e-3 & 1.4289e-4 & 2473.4\\
    \bottomrule
  \end{tabular}
\end{table}

\newcommand\scale[5][0.1]{\parbox{#2\linewidth}{\vspace{0.3pc}%
    \makebox[\linewidth]{\num{#3}\hfill$#4$\hfill\num{#5}}\\[-#1pc]%
    \image[trim=225 1180 275 945,clip=true,width=\linewidth]
    {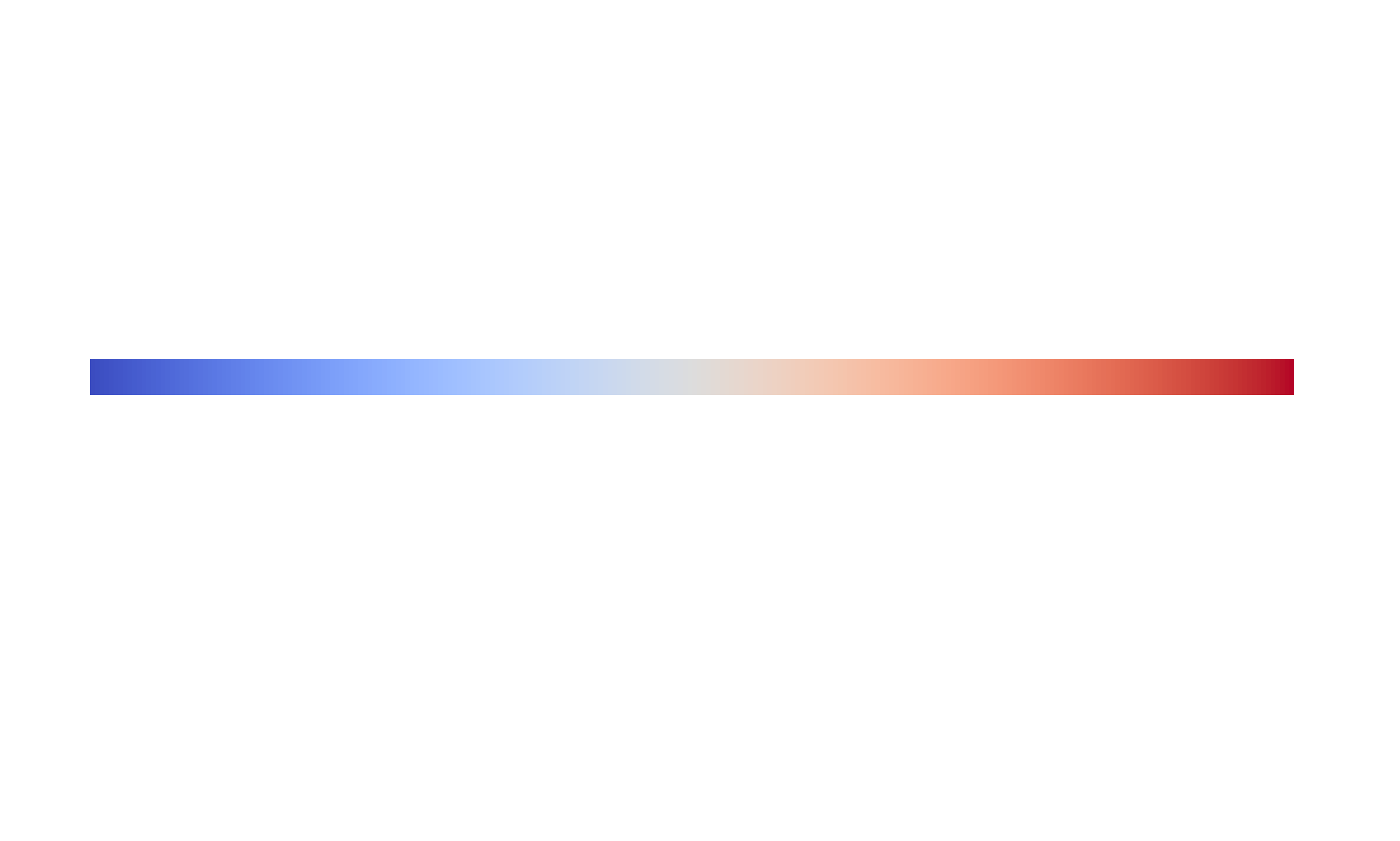}}%
}

\begin{figure}[p]
  \label{fig:phase1}
  \footnotesize
  \centering
  \image[trim=800 175 985 125,clip=true,width=.32\linewidth]
  {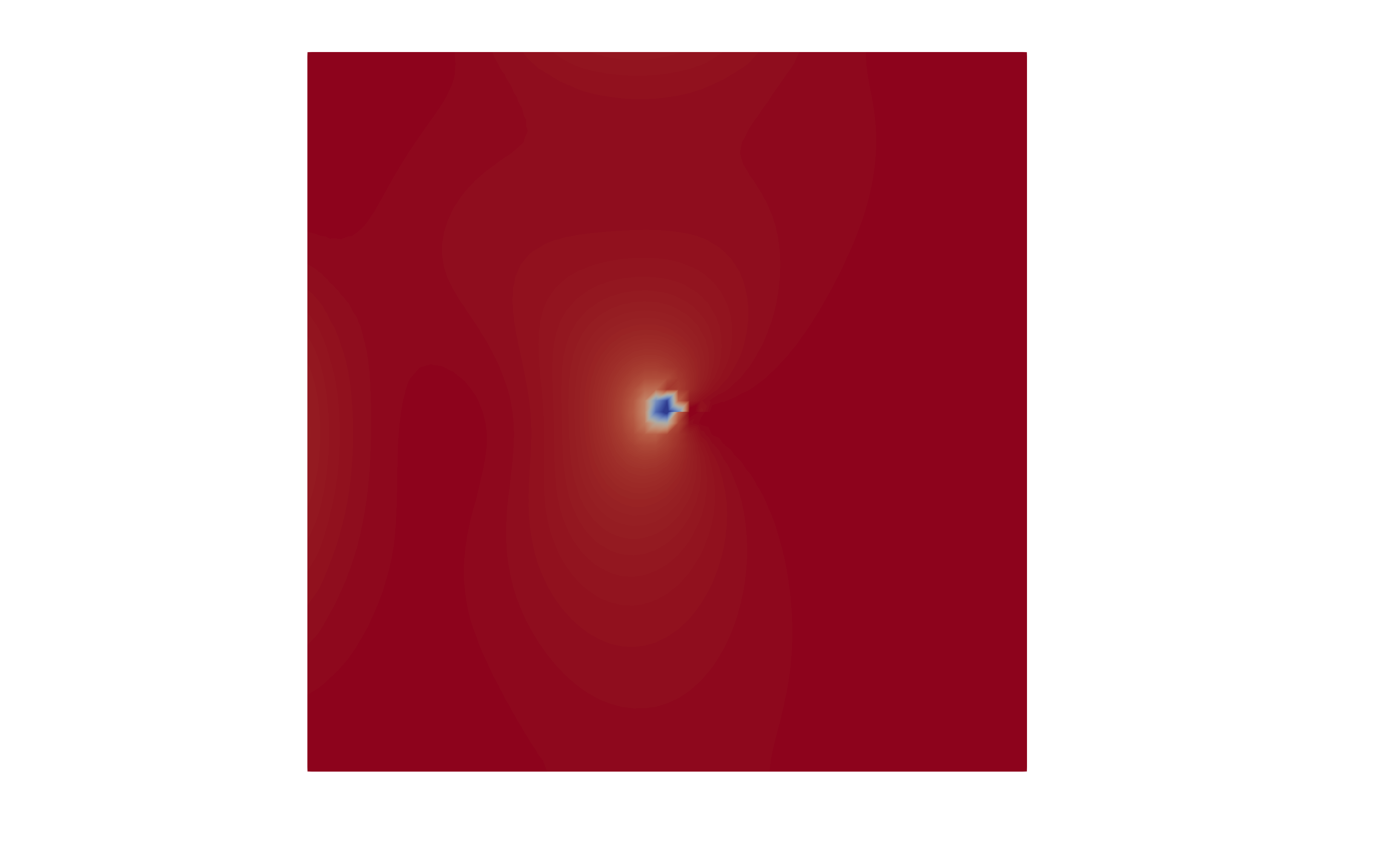}\hfill
  \image[trim=800 175 985 125,clip=true,width=.32\linewidth]
  {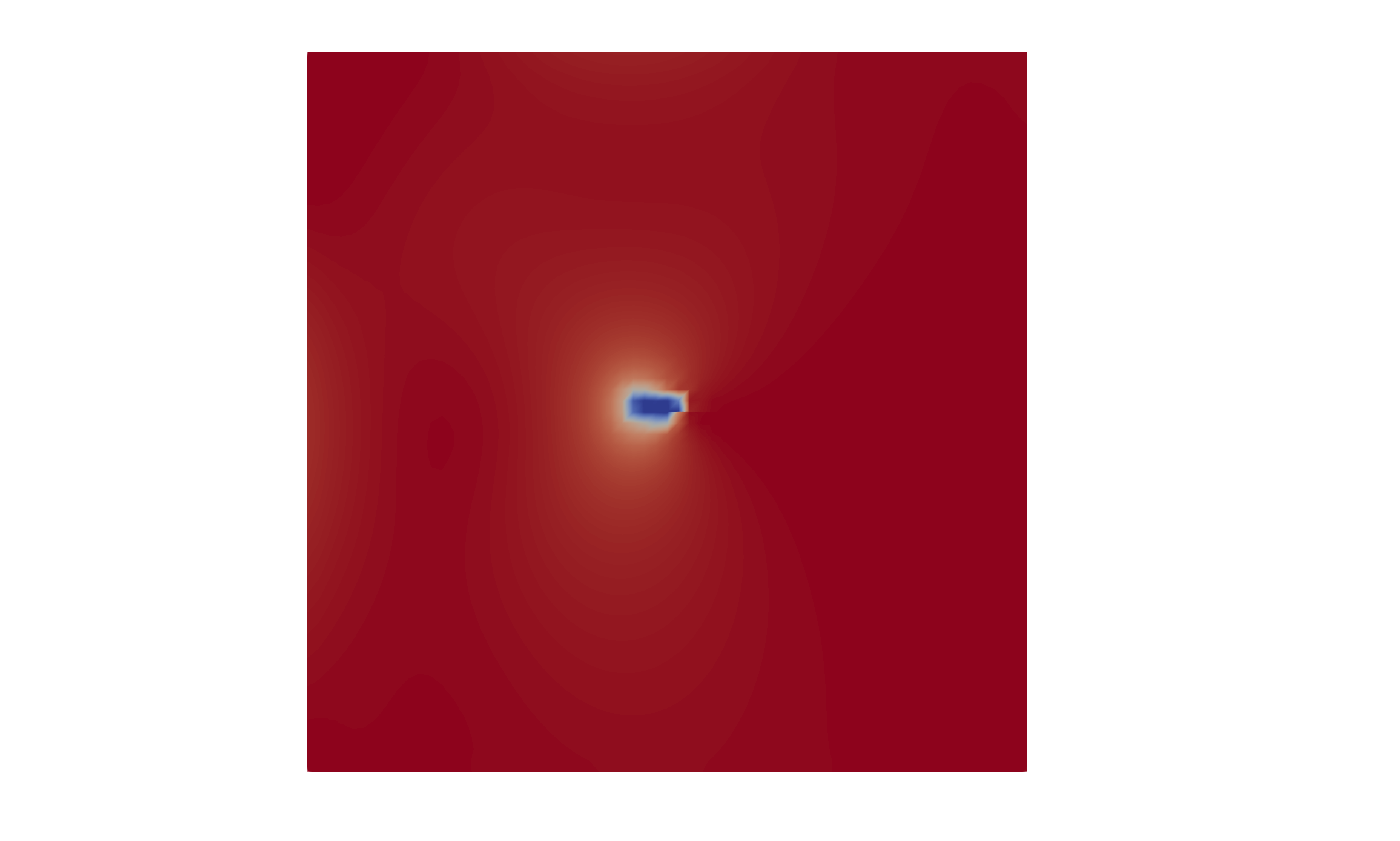}\hfill
  \image[trim=800 175 985 125,clip=true,width=.32\linewidth]
  {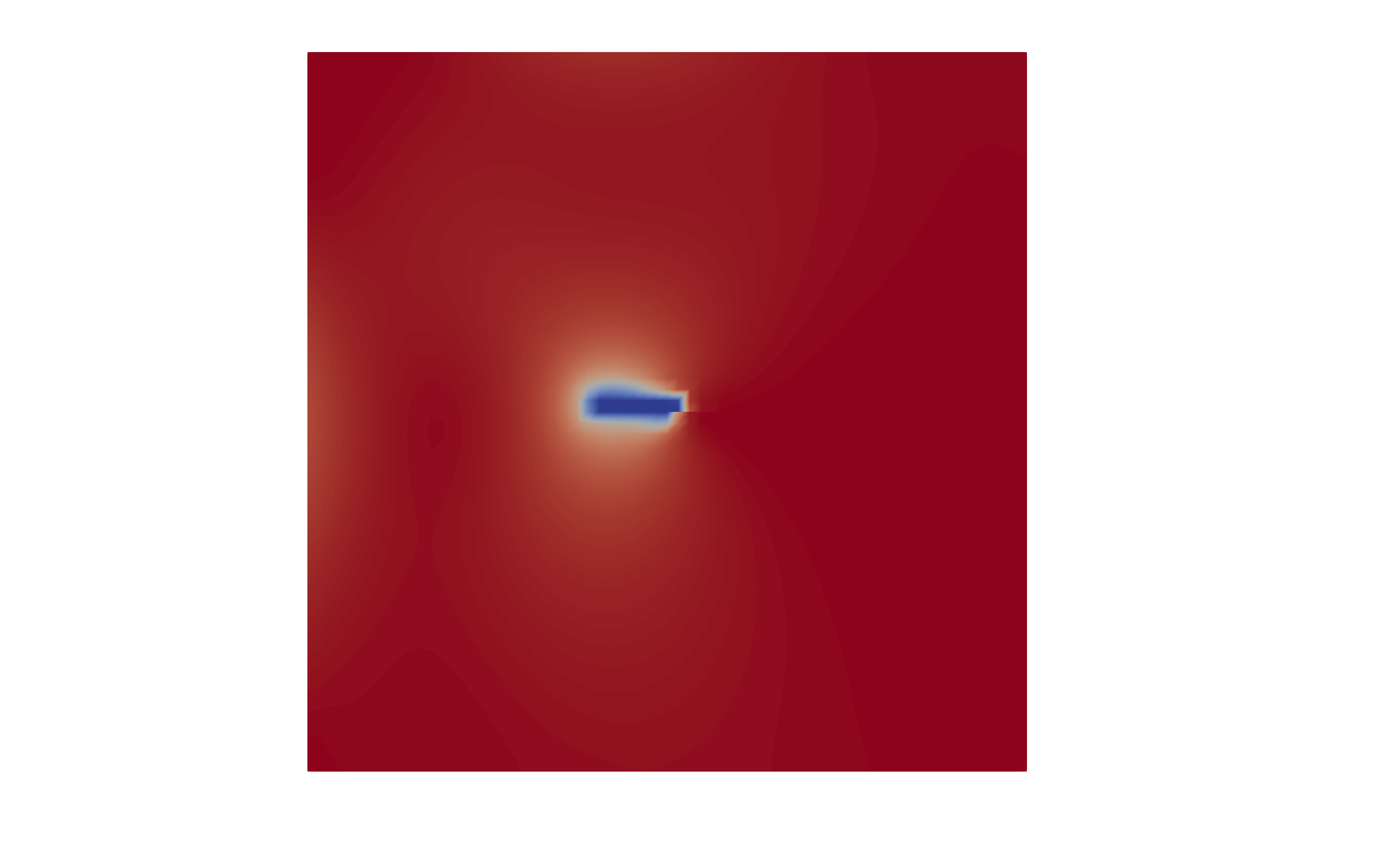}\\%
  \scale[0.3]{0.32}{0}{0.5}{1}%
  \caption{Example 1: optimal phase-field $\varphi$ at times 20, 30, and 40.}
\end{figure}

\begin{figure}[p]
  \label{fig:displacement1}
  \label{fig:adjoint1}
  \footnotesize
  \centering
  \image[trim=1062 325 1070 320,clip=true,width=0.47\linewidth]
  {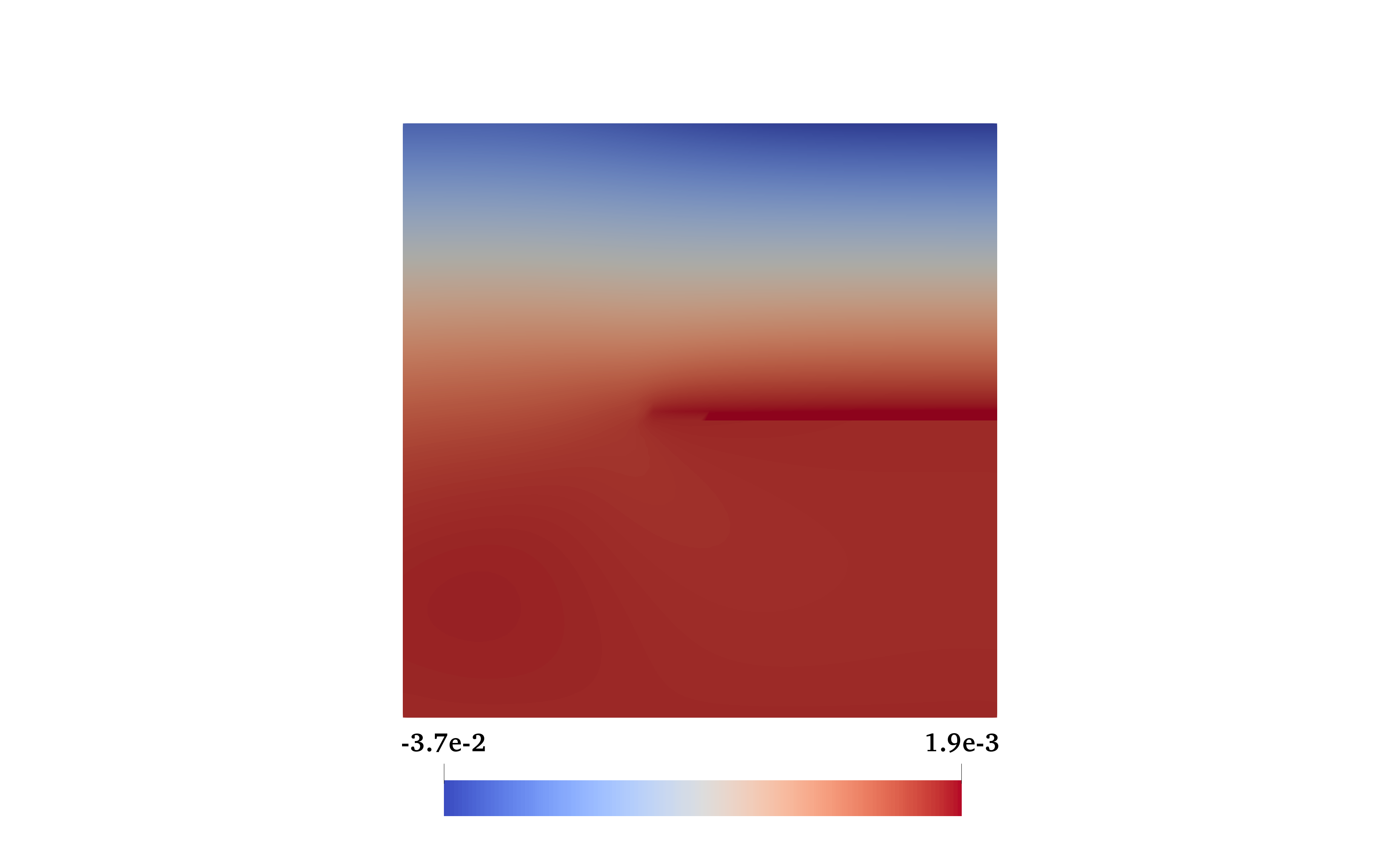}\hfill
  \image[trim=1062 325 1070 320,clip=true,width=0.47\linewidth]
  {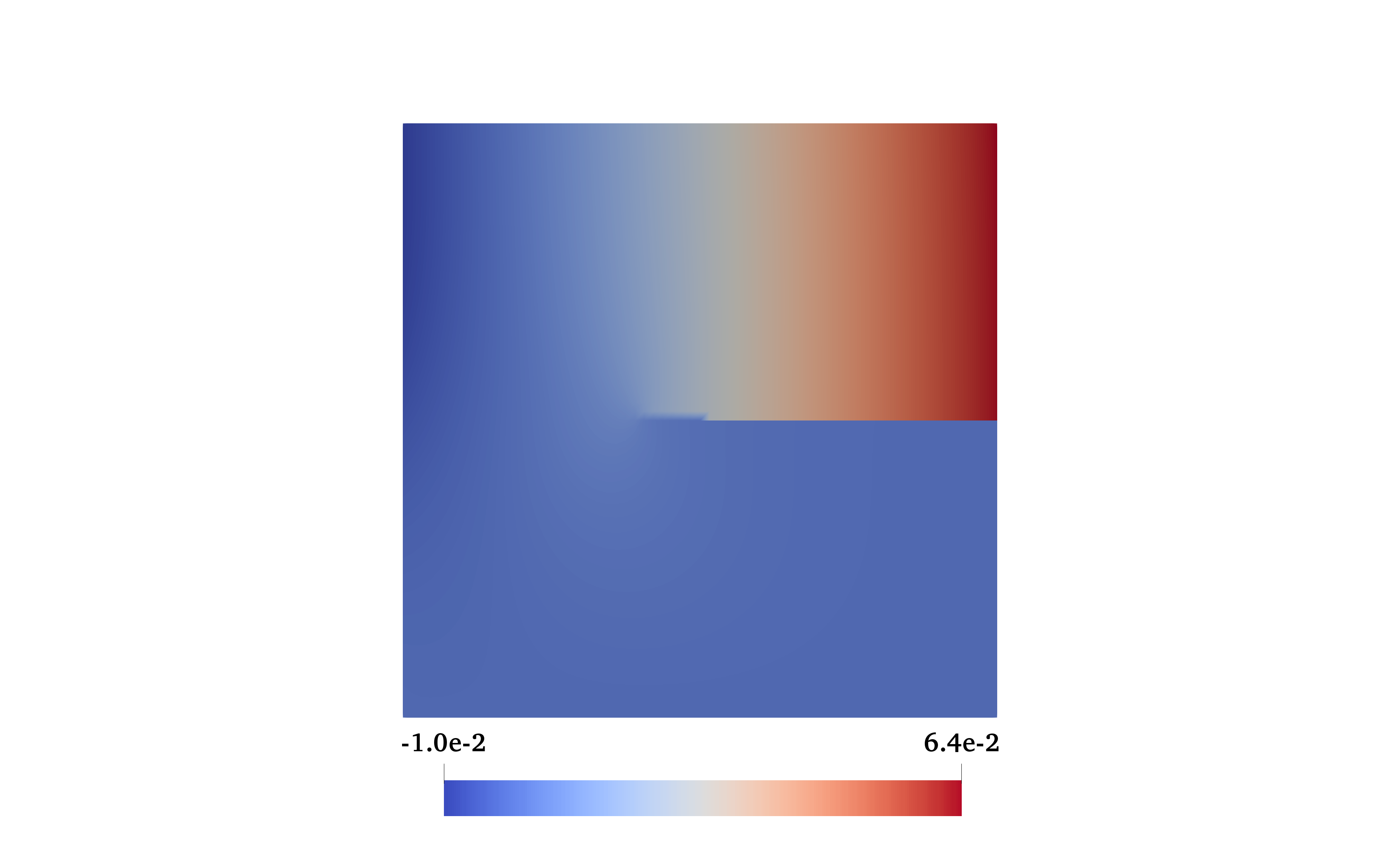}\\%
  \scale{0.47}{-3.7e-2}{}{1.9e-3}\hfill
  \scale{0.47}{-1.0e-2}{}{6.4e-2}\\[3ex]%
  \image[trim=1062 325 1070 320,clip=true,width=0.47\linewidth]
  {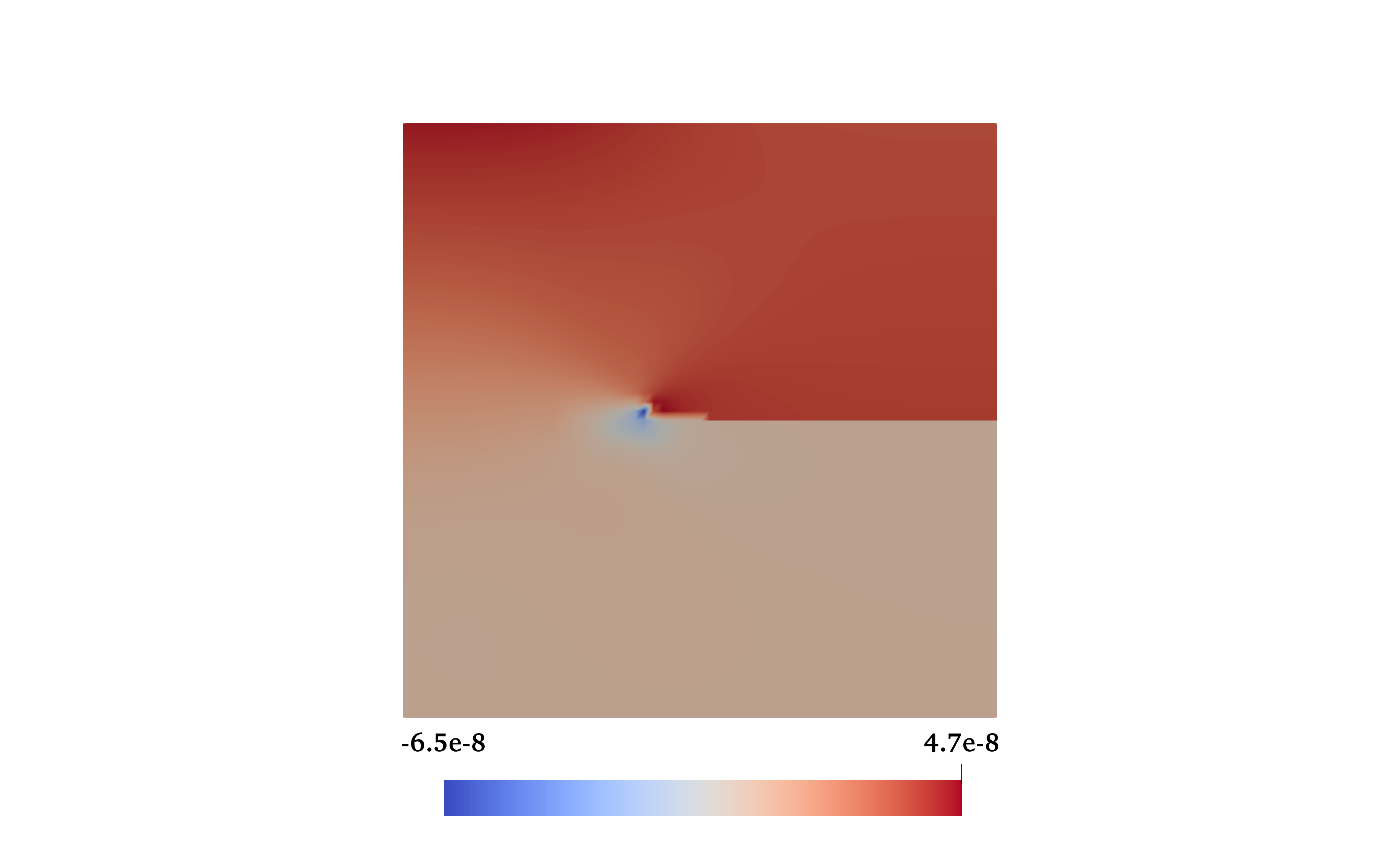}\hfill
  \image[trim=1062 325 1070 320,clip=true,width=0.47\linewidth]
  {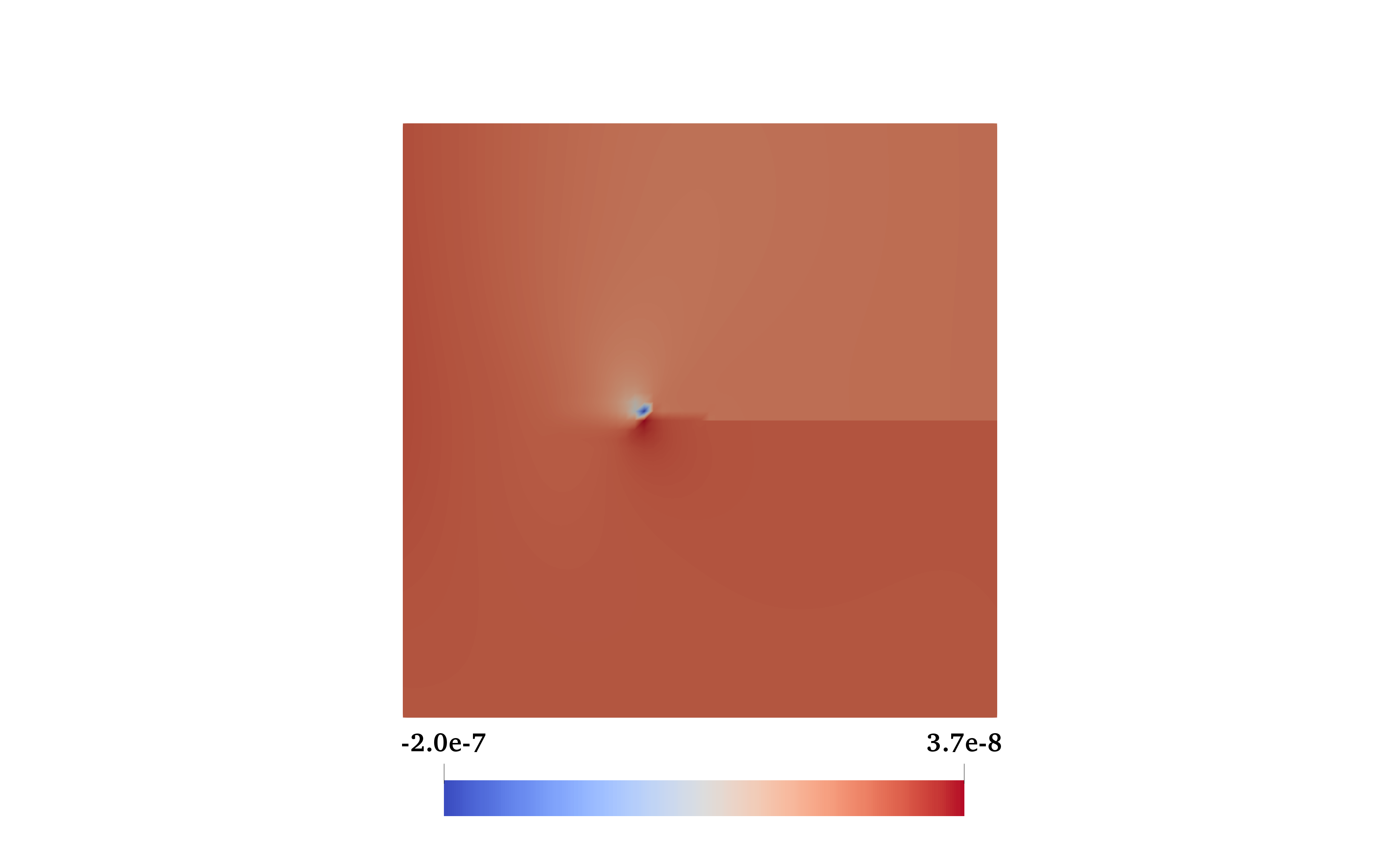}\\
  \scale{0.47}{-6.5e-8}{}{4.7e-8}\hfill
  \scale{0.47}{-2.0e-7}{}{3.7e-8}%
  \caption{Example 1: optimal displacement field $u$ (top: $x$ left, $y$ right)
    and adjoint field $z_u$ (bottom: $x$ left, $y$ right) at time 40.}
\end{figure}

\begin{table}
  \label{tab:results2}
  \centering
  \caption{Results of Example 2: numerical solver performance,
      convergence of functionals, and evolution of control forces.}
  \sisetup{table-format=1.4e+1}
  \begin{tabular}{rr
    S[table-format=1.2e+1]S[table-format=1.2e+2]SSSS[table-format=4.1]}
    \toprule
    Iter&CG&{Relative}&{Absolute}&{Cost}&{Tracking}&{Tikhonov}&{Force}\\[-2pt]
        &  &{residual}&{residual}&      &          &       &          \\
    \midrule
    0 & ---& 1.0     & 2.92e-07 & 1.1481e-2 & 1.1037e-2 & 4.4409e-4 & 1.0\\
    1 &  3 & 0.475   & 1.39e-07 & 1.0703e-2 & 1.0703e-2 & 2.0512e-9 & 3001.0\\
    2 &  3 & 0.280   & 8.17e-08 & 1.0308e-2 & 1.0208e-2 & 1.0027e-4 & 5542.4\\
    3 &  3 & 0.011   & 3.21e-09 & 1.0068e-2 & 9.8146e-3 & 2.5294e-4 & 7051.8\\
    4 & 14 & 2.02e-3 & 5.89e-10 & 9.9970e-3 & 9.6262e-3 & 3.4301e-4 & 7698.9\\
    5 & 11 & 3.59e-4 & 1.05e-10 & 9.9953e-3 & 9.5933e-3 & 3.5995e-4 & 7767.2\\
    6 &  9 & 2.14e-4 & 6.26e-11 & 9.9951e-3 & 9.5880e-3 & 3.6268e-4 & 7788.6\\
    7 &  6 & 1.60e-4 & 4.66e-11 & 9.9949e-3 & 9.5855e-3 & 3.6402e-4 & 7791.8\\
    8 &  8 & 1.07e-4 & 3.12e-11 & 9.9948e-3 & 9.5827e-3 & 3.6542e-4 & 7802.2\\
    9 &  5 & 6.93e-5 & 2.02e-11 & 9.9947e-3 & 9.5811e-3 & 3.6625e-4 & 7805.3\\
    10 & 6 & 4.32e-5 & 1.26e-11 & 9.9947e-3 & 9.5799e-3 & 3.6688e-4 & 7809.6\\
    11 & 4 & 2.99e-5 & 8.73e-12 & 9.9946e-3 & 9.5792e-3 & 3.6723e-4 & 7810.9\\
    \bottomrule
  \end{tabular}
\end{table}

As a side-note, we observed within our computations that
it is quite challenging to entirely reach the desired crack $\varphi_d$.
This is clear since we have the usual competition
between the physics given by the tracking functional
and the numerical regularization given by the Tikhonov term.
Of course, with $\alpha$ chosen sufficiently well
and a well-guessed $q_d$ we could place a higher weight
on the regularization to track the desired fracture path better.
On the other hand, to control a propagating fracture at all
has not yet been achieved in the published literature
to the best of our knowledge.

\subsection{Example 2: horizontal fracture in the middle}
The second example is motivated by the question whether
it is possible to produce a one sided crack growth.
The desired phase-field $\varphi_d$ continues the initial notch
only to the left boundary, see again \cref{fig:examples},
but the right fracture tip should not move.
For that reason, this numerical experiment differs
quite significantly from the first configuration in the behavior
of the numerical solution as well as the final fracture path outcomes.

\begin{figure}[p]
  \label{fig:phase2}
  \footnotesize
  \centering
  \image[trim=895 155 905 150,clip=true,width=.32\linewidth]
  {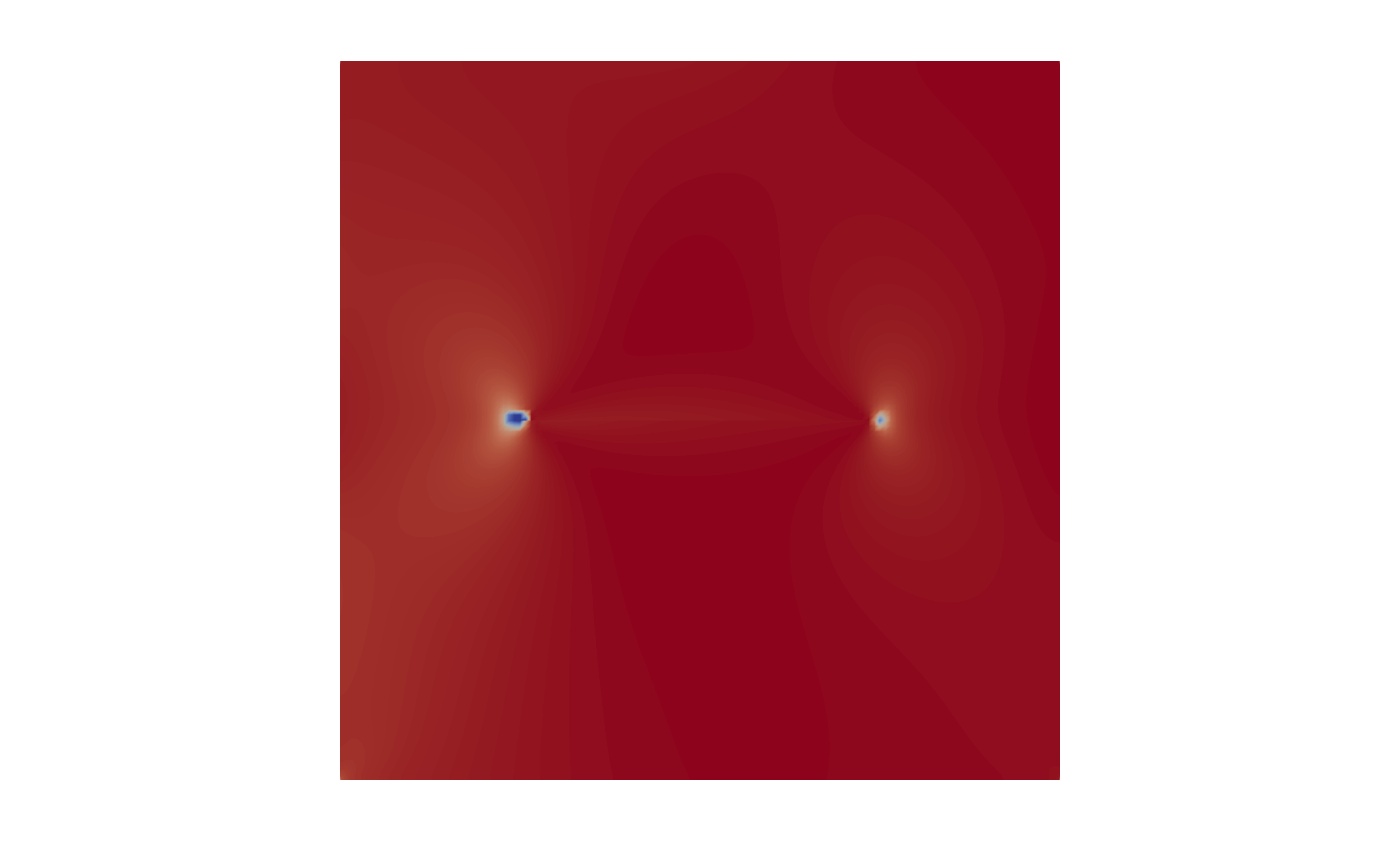}\hfill
  \image[trim=895 155 905 150,clip=true,width=.32\linewidth]
  {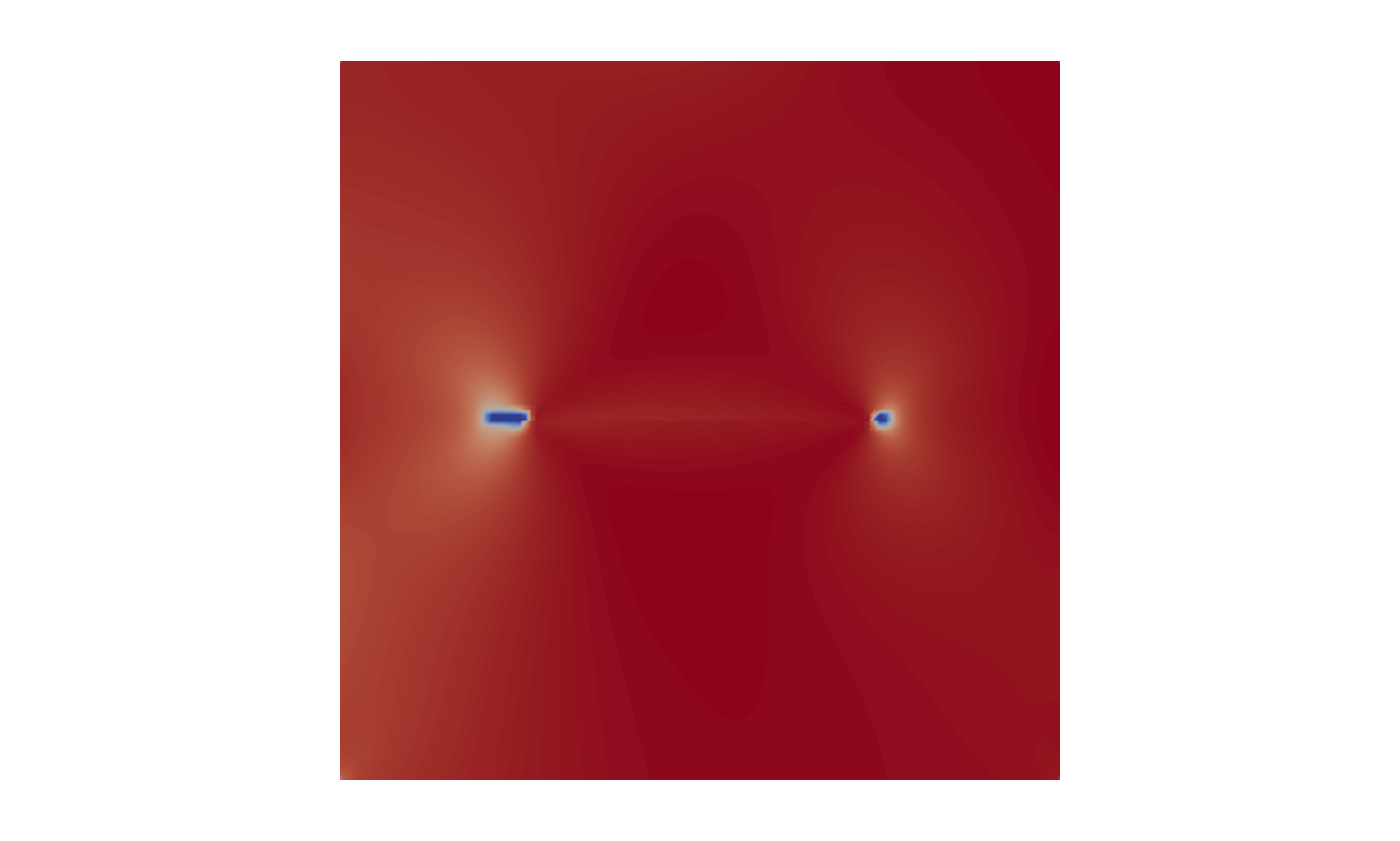}\hfill
  \image[trim=895 155 905 150,clip=true,width=.32\linewidth]
  {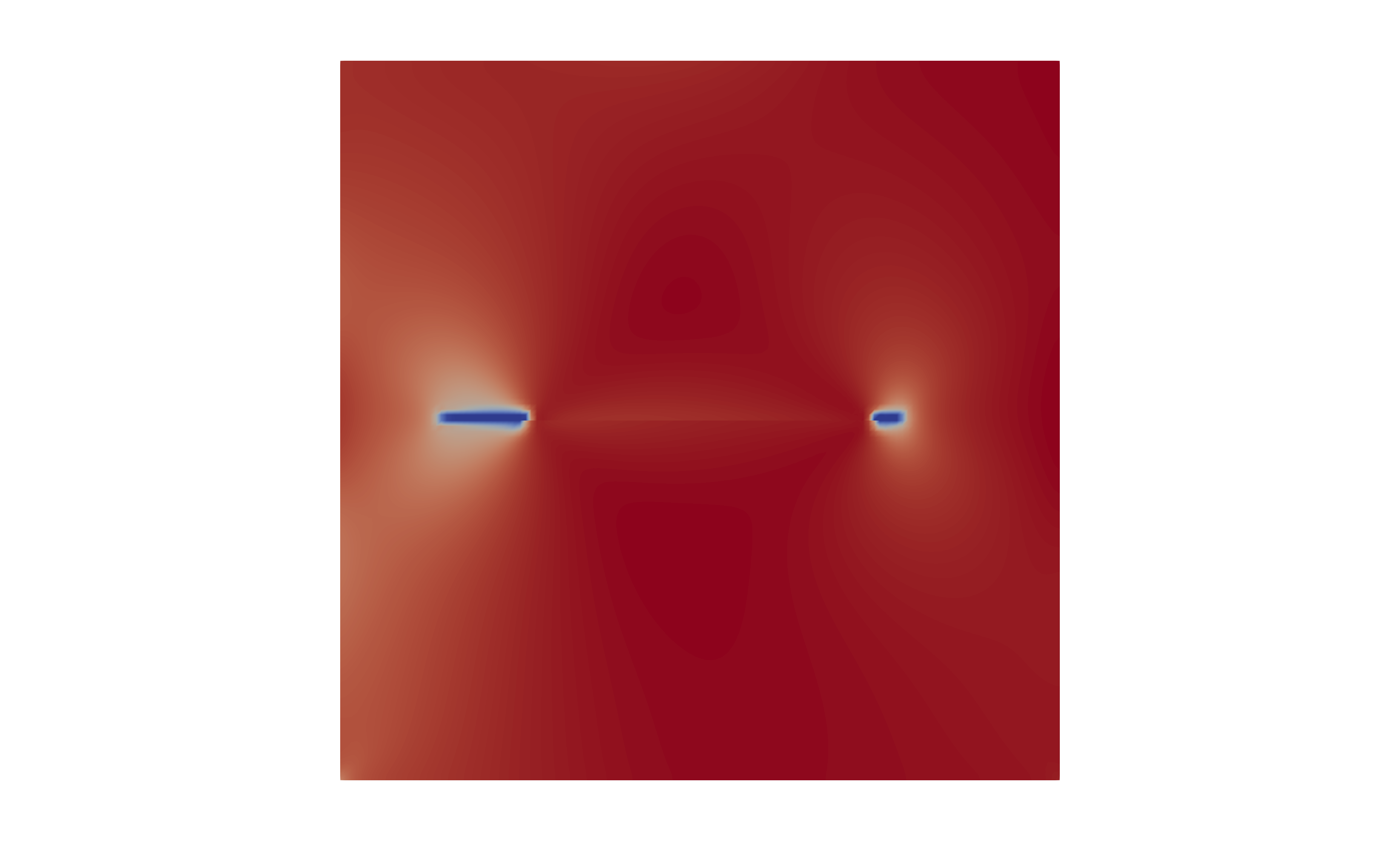}\\%
  \scale[0.3]{0.32}{0}{0.5}{1}%
  \caption{Example 2: optimal phase-field $\varphi$ at times 20, 30, and 40.}
\end{figure}

\begin{figure}[p]
  \label{fig:displacement2}
  \label{fig:adjoint2}
  \footnotesize
  \centering
  \image[trim=1005 325 1010 325,clip=true,width=0.47\linewidth]
  {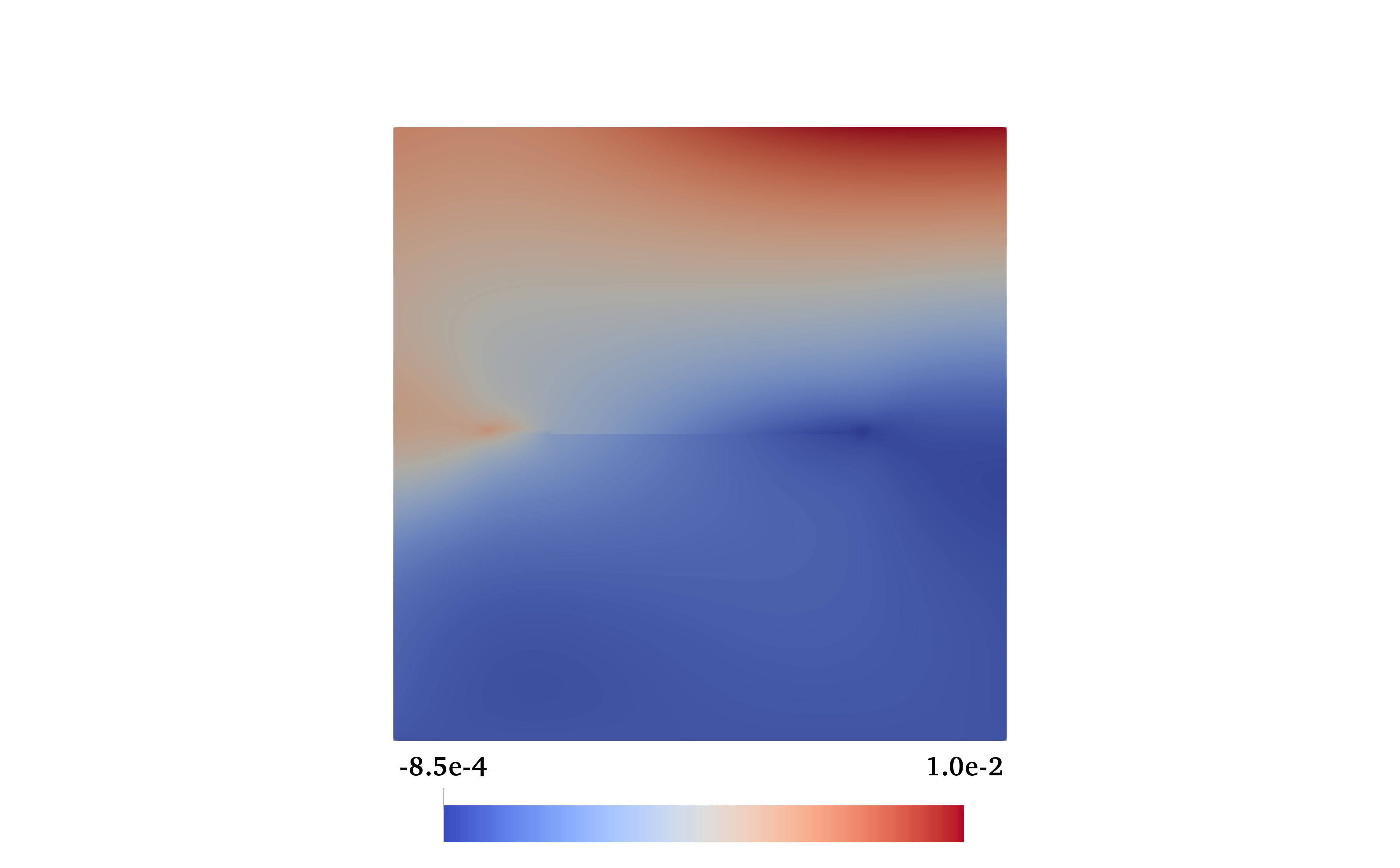}\hfill
  \image[trim=1005 325 1010 325,clip=true,width=0.47\linewidth]
  {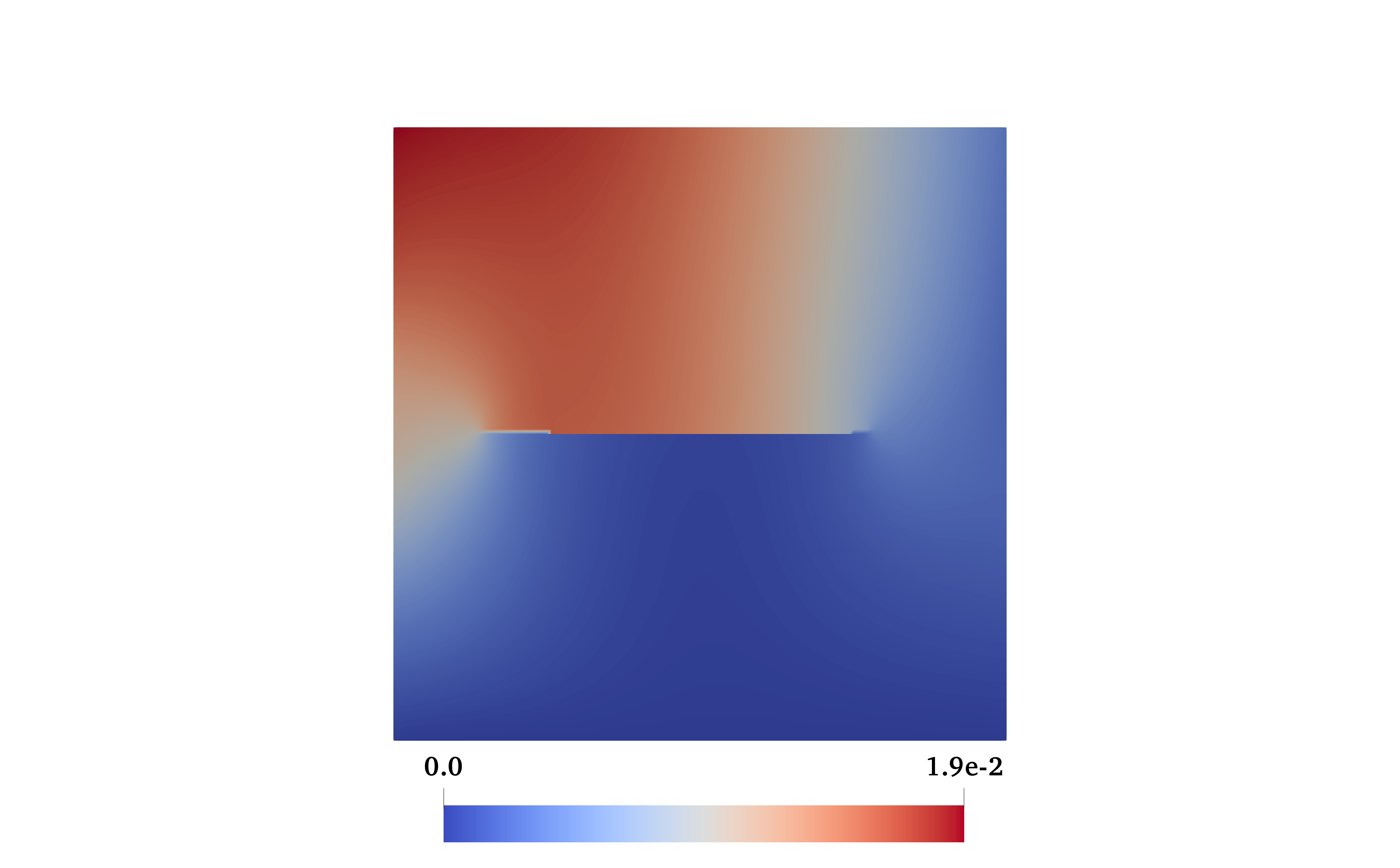}\\%
  \scale{0.47}{-8.5e-4}{}{1.0e-2}\hfill
  \scale{0.47}{0.0}{}{1.9e-2}\\[3ex]%
  \image[trim=1005 325 1010 325,clip=true,width=0.47\linewidth]
  {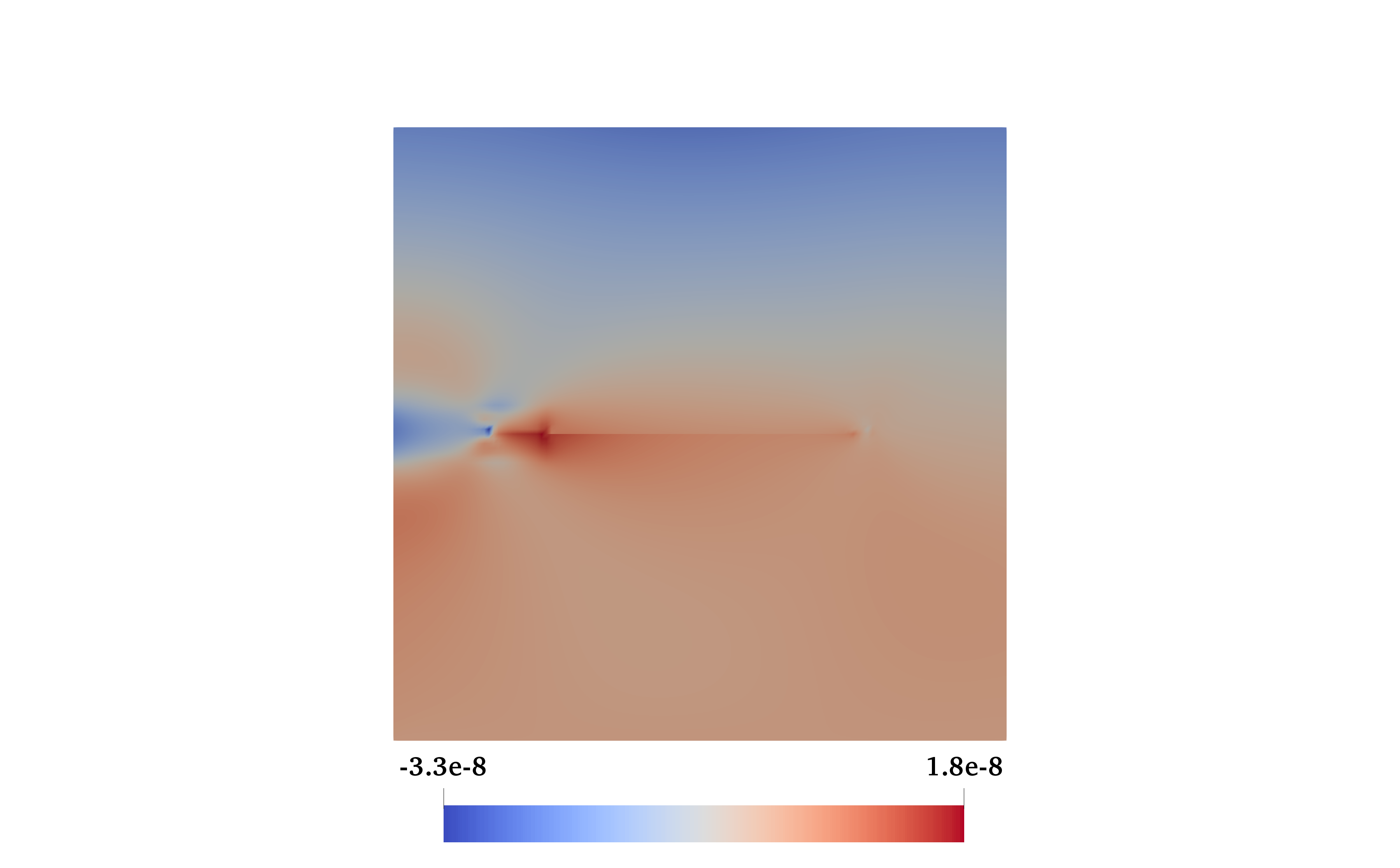}\hfill
  \image[trim=1005 325 1010 325,clip=true,width=0.47\linewidth]
  {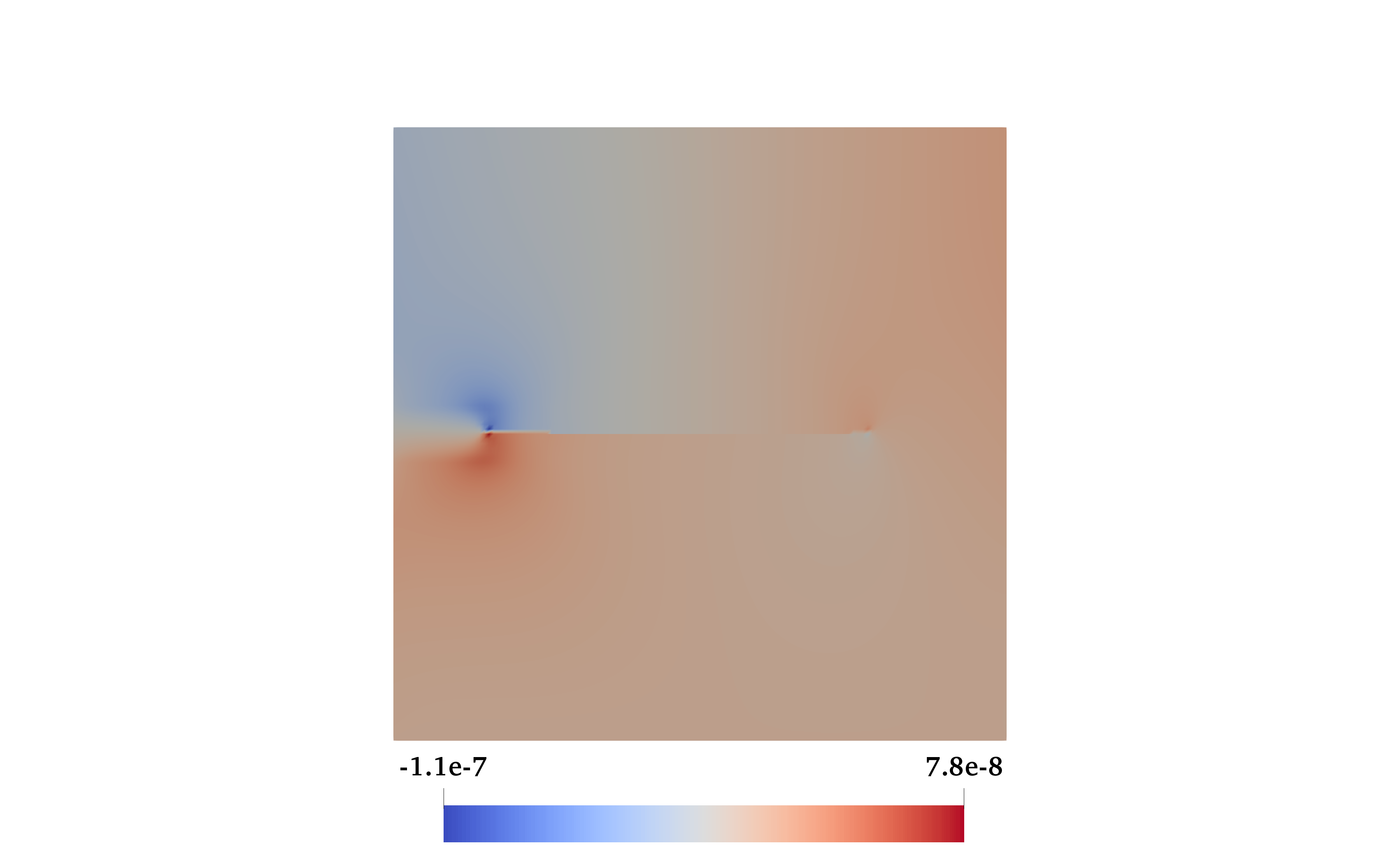}\\%
  \scale{0.47}{-3.3e-8}{}{1.8e-8}\hfill
  \scale{0.47}{-1.1e-7}{}{7.8e-8}%
  \caption{Example 2: optimal displacement field $u$ (top: $x$ left, $y$ right)
    and adjoint field $z_u$ (bottom: $x$ left, $y$ right) at time 40.}
\end{figure}

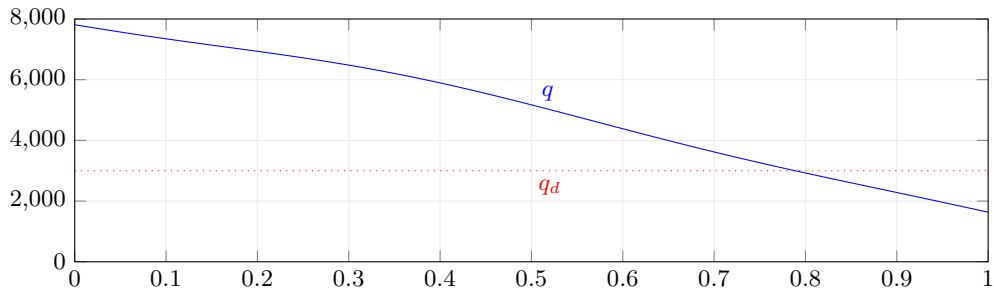
\begin{figure}
  \label{fig:force2}
  \footnotesize
  \centering
  \begin{tikzpicture}
    \begin{axis}[width=0.8\linewidth,height=48mm,
      xmin=0,xmax=1,ymin=0,ymax=8000,ytick distance={2000},grid=major]
      \draw[red,dotted] (0,3000) -- (1,3000) node[pos=0.52,below] {$q_d$};
      \addplot[blue] table {Figures/Test_2/Kraft.txt} node[pos=0.45,above] {$q$};
    \end{axis}
  \end{tikzpicture}
  \caption{Example 2: optimal control force (solid) and nominal control
    force (dotted) on upper boundary $\Gamma_N = [0, 1] \times \set{1}$.}
\end{figure}

Our results are presented in
\cref{tab:results2,fig:force2,fig:phase2,fig:displacement2,fig:adjoint2}.
Here the tolerance for Newton's method is \num{1e-11}.
The qualitative performance of each numerical solver
is similar to Example 1.
Concerning the optimal fracture path (phase-field in \cref{fig:phase2})
we obtain the desired solution,
namely crack growth starting from the left fracture tip.
In the rightmost subfigure of \cref{fig:phase2}, however,
we also observe a slight movement of the right fracture tip.
This is reasonable because the optimal traction
still has a physical impact on the overall fracture and the cost functional
is only enforced up to numerical regularization and discretization
approximation qualities.
The optimal force now decreases from 7810.9 at $(0, 1)$ to 1633.6 at $(1, 1)$.

\section{Conclusions}
\label{sec_conclusions}
In this paper we derived a space-time Galerkin
formulation for a regularized phase-field fracture optimal control setting.
By introducing jump terms in time and with the help of a discontinuous
Galerkin discretization in time, specific time-stepping schemes could
be obtained. A careful investigation of correct weighting
of two regularization terms and the initial conditions was necessary
for the forward phase-field fracture problem. The solution process
of the optimization problem was based on the reduced approach
in which the state variables are obtained from a solution operator
acting on the controls. The numerical solution algorithm is based
on Newton's method in which three auxiliary problems are required.
The main part of the paper was concerned with the detailed derivation
of these terms, which are to the best of our knowledge novel
in the published literature.
We then discussed two numerical tests in order to show the performance
of our framework. Therein, we studied the convergence of the
residuals (relative and absolute) of the reduced problem,
the number of Newton steps and CG iterations as well as the
behavior of the cost functional. These findings indicate the
robustness and suitability of our numerical solvers to address
optimal control phase-field with propagating fractures.
Graphical results of the phase-field solution (showing the
fracture path) and the optimal displacement field illustrate
our findings.

\section{Acknowledgements}
The first and third author are partially funded by the
      Deutsche Forschungsgemeinschaft (DFG, German Research Foundation)
      Priority Program 1962 (DFG SPP 1962) within the subproject
      \textit{Optimizing Fracture Propagation using a Phase-Field Approach}
      with the project number 314067056.
      The second author is funded by the DFG -- SFB1463 -- 434502799.


\end{document}